\documentclass[12pt]{article}
\usepackage{amsmath}
\usepackage{pifont}
\usepackage{mathrsfs}
\usepackage{amsfonts}
\usepackage{stmaryrd}
\usepackage{latexsym,amssymb,amsmath}
\usepackage[colorlinks,linkcolor=black,anchorcolor=black,citecolor=black]{hyperref}
\newcommand{\mbb}{\mathbb}
\newcommand{\ol}{\overline}
\newcommand{\psp}{\vspace{0.2cm}}
\newcommand{\td}{\tilde}

\begin{document}

\baselineskip=18pt
\newcommand{\headingstobeshown}{}
\def\tenrm{\rm}

\renewcommand{\theequation}{\thesection.\arabic{equation}}

%\makeatletter      % '@' is now a normal "letter" for TeX
%\@addtoreset{equation}{section}
%\makeatother       % '@' is restored as a "non-letter" character for TeX

\newenvironment{proof}
               {\begin{sloppypar} \noindent{\it Proof.}}
               {\hspace*{\fill} $\square$ \end{sloppypar}}

\newtheorem{atheorem}{\bf \temp}[section]

\newenvironment{theorem}[1]{\def \temp{#1}
                \begin{atheorem}\medskip
                }
                {\medskip
                 \end{atheorem}}

\newtheorem{thm}[atheorem]{Theorem}
\newtheorem{cor}[atheorem]{Corollary}
\newtheorem{lem}[atheorem]{Lemma}
\newtheorem{pro}[atheorem]{Property}
\newtheorem{prop}[atheorem]{Proposition}
\newtheorem{de}[atheorem]{Definition}
\newtheorem{rem}[atheorem]{Remark}
\newtheorem{fac}[atheorem]{Fact}
\newtheorem{ex}[atheorem]{Example}
\newtheorem{pr}[atheorem]{Problem}
\newtheorem{cla}[atheorem]{Assert}
\newcommand{\ptl}{\partial}
\newcommand{\lmd}{\lambda}
\newcommand{\la}{\langle}
\newcommand{\ra}{\rangle}

\begin{center}{\Large \bf On Polynomial Representations of Classical}\end{center}
\begin{center}{\Large \bf Strange Lie Superalgebras}\end{center}
\vspace{0.2cm}

\begin{center}{\large Cuiling Luo}
\end{center}

\begin{center}{Institute of Mathematics, Academy of Mathematics \&
System Sciences,
%}\end{center}
%\begin{center}{
\\Chinese Academy of Sciences, Beijing 100190, China}\end{center}
\begin{center}{E-mail: luocuiling@amss.ac.cn}\end{center}
\vspace{0.4cm}
\begin{abstract}

 In this paper, various polynomial representations of
strange classical Lie superalgebras are investigated. It turns out
that the representations for the algebras of type $P$ are
indecomposable, and we obtain the composition series of the
underlying modules. As  modules of the algebras of type $Q$ , the
polynomial algebras are decomposed into a direct sum of
irreducible submodules.

\vspace{0.5cm} {\bf Keywords:} polynomial, representation, Lie
superalgebra, composition series.
\end{abstract}
\vspace{0.5cm}

\section{Introduction}

Lie superalgebras were introduced by physicists as the fundamental
tools of studying the supersymmetry in physics. Unlike Lie algebra
case, finite-dimensional modules of finite-dimensional simple Lie
superalgebras may not be completely reducible and the structure of
finite-dimensional irreducible modules is much more complicated
due to the existence of so-called {\it atypical} modules (cf.
[17], [18]). In his celebrated work [16] on classification of
finite-dimensional simple Lie superalgebras, Kac found two
families of exotic classical simple Lie superalgebras, which are
called ``strange" Lie superalgebras of type $P$ and $Q$,
respectively. These superalgebras do not have analogues in Lie
algebras. The strange Lie superalgebras have attracted a number of
mathematicians' attention.

Javis and  Murray [9] obtained the Casimir invariants,
characteristic identities, and tensor operators for the strange
Lie superalgebras. Moreover, Nazarov [11] found Yangians of the
superalgebras. In [2], Frappat, Sciarrino and  Sorba studied
Dynkin-like diagrams and a certain representation of the strange
superalgebra $P(n)$. In addition, they [3] gave a certain
oscillator realization of the strange superalgebras. Medak [13]
generalized the Baker-Campbell-Hausdorff formula and used it to
examine the so-called BCH-Lie and BCH-invertible subalgebras in
the Lie superalgebra $P(n)$. Penkov and Serganova [20] discovered
a surprising phenomena that the multiplicity of the highest weight
in the finite-dimensional irreducible representations of $q(n)$ is
in general greater than 1.

Gruson  [8] computed the cohomology with trivial coefficients for
the strange Lie superalgebras. Palev and  Van der Jeugt [19] found
a family of nongraded  Fock representations of the Lie
superalgebra $Q(n)$. Gorelik [5] obtained the center of the
universal enveloping algebra of the Lie superalgebra of type $P$.
Serganova [21] determined the center of the quotient algebra of
the universal enveloping algebra of the Lie superalgebra of type
$P$ by its Jacobson radical and used it to study the typical
highest weight modules of the algebra. Medak [14] proved that each
maximal invertible subalgebra of $P(n)$ is $\mathbb{Z}$-graded.
Moon [15] obtained a ``Schur-Weyl duality" for the algebras of
type $P$. Martinez and Zelmanov [12] classified Lie superalgebras
graded by $P(n)$ and $Q(n)$.  Brundan [1] found a connection
between Kazhdan-Lusztig polynomials and character formulas for the
Lie superalgebra $q(n)$. Gorelik [6] obtained the Shapovalov
determinants of Q-type Lie superalgebras. Stukopin [22] studied
the Yangians of the strange Lie superalgebra of type $Q$ by
Drinfel'd approach. Gorelik and Serganova [7] investigated the
structure of Verma modules over the twisted affine Lie
superalgebra $q(n)^{(2)}$.

One way of understanding simple Lie algebras and simple Lie
superalgebras is to determine the structure of their natural
representations. Canonical polynomial irreducible representations
(also known as oscillator representations in physics (e.g., cf.
[4])) of finite-dimensional simple Lie algebras are very important
from application point of view, where both the representation
formulas and bases are clear. In [10], we determined the structure
of certain noncanonical polynomial representations of classical
simple Lie algebras, in particular, their irreducible submodules.
In this paper, we want to generalize the above results to the
strange simple Lie superalgebras. The details are as follows.

Throughout this paper, we denote by $\mbb{Z}$ the ring of integers
and by $\mbb{N}$ the set of nonnegative integers. For convenience,
we also use the following notation of indices:
\begin{equation}
\ol{i,j}=\{i,i+1,...,j\},\end{equation} where $i\leq j$ are
integers. Let $E_{i,j}$ be the matrix whose $(i,j)$-entry is 1 and
the others are zero. Moreover, all the  vector spaces are assumed
over $\mathbb{C}$, the field of complex numbers. The general
linear Lie superalgebra $gl(m|n)=gl(m|n)_{\bar{0}}\oplus
gl(m|n)_{\bar{1}}$ with
\begin{equation}gl(m|n)_{\bar{0}}=\sum_{i,j=1}^m\mathbb{C}E_{i,j}+\sum_{p,q=1}^n
\mathbb{C}E_{m+p,m+q},\;
gl(m|n)_{\bar{1}}=\sum_{i=1}^m\sum_{p=1}^n(\mathbb{C}E_{i,m+p}+
\mathbb{C}E_{m+p,i}).\end{equation} Let ${\cal
A}=\mathbb{C}[x_1,\cdots,x_m,\theta_1,\cdots,\theta_n]$ be the
polynomial algebra in bosonic variables $x_1,\cdots,x_m$ and
fermionic variables $\theta_1,\cdots,\theta_n$, i.e.
\begin{equation}
x_ix_j=x_jx_i,\;\;\theta_p\theta_q=-\theta_q\theta_p,\;\;x_i\theta_p=\theta_px_i,\;\;i,j\in\overline{1,m},\;\;p,q\in\overline{1,n}.
\end{equation}

Take an integer $r\in\overline{0,m}$. Define an action of $gl(m|n)$
on ${\cal A}$ by
\begin{equation}
\begin{array}{l}
E_{i,j}|_{\mathcal{A}}=\left\{\begin{array}{lll}
-x_j\partial_{x_i}-\delta_{i,j} &\mbox{if}&i,j\in\ol{1,r},\\
\partial_{x_i}\partial_{x_j} &\mbox{if} &i\in\ol{1,r},\;j\in\ol{r+1,m},\\
-x_ix_j &\mbox{if} &j\in\ol{1,r},\;i\in\ol{r+1,m},\\
x_i\partial_{x_j}& \mbox{if} &i,j\in\ol{r+1, m},
\end{array}\right.\\ \\
E_{i,m+p}|_{\mathcal{A}}=\left\{
\begin{array}{lll}
\partial_{x_i}\partial_{\theta_p} &\mbox{if} &i\in\ol{1, r},\;p\in\ol{1,n},\\
x_i\partial_{\theta_p} &\mbox{if} &i\in\ol{r+1,m},\;p\in\ol{1,n},
\end{array}\right.\\ \\
E_{m+p,j}|_{\mathcal{A}}=\left\{
\begin{array}{lll}
-\theta_px_j &\mbox{if} &j\in\ol{1, r},\;p\in\ol{1,n},\\
\theta_p\partial_{x_j} &\mbox{if}&j\in\ol{r+1,m},\;p\in\ol{1,n} ,
\end{array}
\right.\\ \\
E_{m+p,m+q}|_{\mathcal{A}}=\theta_p\partial_{\theta_q}\qquad\mbox{for
}\;p,q\in\ol{1,n}.\label{I1}
\end{array}
\end{equation}
Then we obtain a representation $\phi_r:gl(m|n)\rightarrow gl({\cal
A})$, which is obtained from the canonical polynomial representation
by swapping $\ptl_{x_i}$ and $-x_i$ for $i\in\ol{1,r}$. In this
case, the degree operator
\begin{equation}{\cal
D}=\phi_r(\sum_{i=1}^{m+n}E_{i,i})+r=-\sum_{i=1}^rx_i\partial_{x_i}+
\sum_{j=r+1}^mx_j\partial_{x_j}+\sum_{p=1}^n\theta_p\partial_{\theta_p}.\label{I2}\end{equation}
Set
\begin{equation}{\cal A}_k=\{f\in{\cal A}\mid {\cal
D}(f)=kf\},\qquad k\in\mathbb{Z}.\end{equation} It is
straightforward to verify that  all ${\cal A}_k$ are irreducible
$gl(m|n)$-submodules and ${\cal A}=\bigoplus_{k\in\mathbb{Z}}{\cal
A}_k$. This is the starting point of this work.

In this paper, we always assume $n\geq 3$. Recall that the simple
strange Lie superalgebra
\begin{eqnarray}P(n-1)&=&\sum_{i=1}^{n-1}
\mathbb{C}(E_{i,i}-E_{i+1,i+1}-E_{n+i,n+i}+E_{n+i+1,n+i+1})\nonumber\\
& &+\sum_{1\leq i<j\leq
n}[\mathbb{C}(E_{i,j}-E_{n+j,n+i})+\mathbb{C}(E_{j,i}-E_{n+i,n+j})\nonumber\\
& &+
\mathbb{C}(E_{n+i,j}-E_{n+j,i})+\mathbb{C}(E_{i,n+j}+E_{j,n+i})]+\sum_{j=1}^n\mathbb{C}E_{j,n+j}\label{I3}\end{eqnarray}
is a subalgebra of $gl(n|n)$. Denote by $\phi_r^P=\phi_r|_{P(n-1)}$
the restricted representation with $m=n$, by ${\cal A}^r$ the
underlying module ${\cal A}$ and by ${\cal A}_k^r$ its submodule
${\cal A}_k$. Furthermore, the notation $\la S\ra$ stands for the
submodule generated by a set $S$. The notation $\hat{f}$ means that
$f$ is missing in an expression.

\psp

{\bf Theorem 1}. {\it We have the following conclusions on the
representation $\phi_r^P$:

(1) If $r=0$, the subspace ${\cal A}^0_1$ is an irreducible
$P(n-1)$-module and the subspace ${\cal A}^0_k$ with $k>1$ has a
composition series
\begin{equation}
{\cal A}^0_k\supset\la x_1^k\ra\supset\{0\} \qquad\mbox{ if }k\neq
n,
\end{equation}
\begin{equation}
{\cal A}^0_n\supset\la \theta_1\cdots\theta_n\ra\supset\la
x_1^n\ra\supset\{0\}.
\end{equation}

(2) When $r=n$, ${\cal A}^{n}_k=0$ if $k>n$, and the submodules
${\cal A}^n_n$ and ${\cal A}^n_{n-1}$ are irreducible. If $0\neq
k\leq n-2$, the submodule ${\cal A}^{n}_k$ has a composition series
\begin{equation}
{\cal A}^{n}_k\supset\la
x_n^{n-k-2}\sum\limits_{i=1}^n(-1)^ix_i\theta_1\cdots\hat{\theta_i}\cdots\theta_n\ra\supset\{0\}.
\end{equation}
Furthermore, ${\cal A}_0^n$ has a composition series
\begin{equation}
{\cal A}^{n}_0\supset\la1\ra\supset\la
x_n^{n-2}\sum\limits_{i=1}^n(-1)^ix_i\theta_1\cdots\hat{\theta_i}\cdots\theta_n\ra\supset\{0\}.
\end{equation}

(3) In the case of $1\leq r\leq n-1$, we have the following
composition series:
\begin{equation}
{\cal A}_k^{r}\supset \la x_{r+1}^{-k}\ra \supset\{0\}\qquad \mbox{
if }\;n-r\neq k\leq 0,
\end{equation}
\begin{equation}
{\cal A}_k^{r}\supset \la x_r^k\ra \supset\{0\}\qquad \mbox{ if
}\;n-r\neq k>0
\end{equation}
and
\begin{equation}
{\cal A}_{n-r}^{r}\supset\la \theta_{r+1}\cdots\theta_n\ra\supset
\la x_r^{n-r}\ra\supset\{0\}.
\end{equation}
 }\psp

Recall that the Lie superalgebra
\begin{eqnarray}\td
Q(n-1)&=&\sum_{i,j=1}^n\mathbb{C}(E_{i,j}+E_{n+j,n+i})+\sum_{i,j\in\ol{1,n};\;i\neq
j} \mathbb{C}(E_{i,n+j}+E_{n+j,i})\nonumber
\\ & &+\sum_{i=1}^{n-1}\mathbb{C}(E_{i,n+i}-E_{i+1,n+i+1}+E_{n+i,i}-E_{n+i+1,i+1})
\end{eqnarray} is a subalgebra of $gl(n|n)$. Denote by $\phi_r^Q=\phi_r|_{\td
Q(n-1)}$ the restricted representation with $m=n$, by ${\cal A}^r$
the underlying module ${\cal A}$ and by ${\cal A}_k^r$ its
submodule ${\cal A}_k$. Let $I_{2n}$ be the $2n\times 2n$ identity
matrix, i.e. $I_{2n}=\sum_{i=1}^{2n}E_{i,i}$. Then
$\mathbb{C}I_{2n}$ is the center of $\td Q(n-1)$ and the simple
strange Lie superalgebra
\begin{equation}Q(n-1)=\td
Q(n-1)/\mathbb{C}I_{2n}.\end{equation}\psp

{\bf Theorem 2}. {\it We have the following conclusions on the
representation $\td\phi_r^Q$ of $\td Q(n-1)$:

(1) For $k>0$, the submodule ${\cal A}^0_k$ is a direct sum of two
irreducible submodules $\la x_1^k+\sqrt{k}x_1^{k-1}\theta_1\ra$
and $\la x_1^k-\sqrt{k}x_1^{k-1}\theta_1\ra$.

(2) The submodules ${\cal A}_k^{r}$ with $k\in\mathbb{Z}$ are
irreducible if $0<r<n$.

(3) The representations $\td\phi_0^Q$ and $\td\phi_n^Q$ are
equivalent}.\psp

Since $\phi_r^Q(I_{2n})|_{{\cal A}^r_r}={\cal D}|_{{\cal
A}^r_r}-r=0$ (cf. (\ref{I2})), the subspaces ${\cal A}^r_r$ with
$0<r<n$ are irreducible submodules of the  strange simple Lie
superalgebra $Q(n-1)$ by (1.16).

The proof of Theorem 1 is given in Section 2. In Section 3, we prove
Theorem 2.

\section{Proof of Theorem 1}
\setcounter{equation}{0}

In this section, we investigate the representations $\phi_r^P$  of
the strange simple Lie superalgebra $P(n-1)$ in (\ref{I3}).

Note that
\begin{equation}H=\sum_{i=1}^{n-1}
\mathbb{C}(E_{i,i}-E_{i+1,i+1}-E_{n+i,n+i}+E_{n+i+1,n+i+1})\end{equation}
forms a Cartan subalgebra of the even subalgebra
\begin{equation}P(n-1)_{\bar{0}}=P(n-1)\cap gl(n|n)_{\bar 0}=
H +\sum_{i,j\in\ol{1,n};\;i\neq
j}\mathbb{C}(E_{i,j}-E_{n+j,n+i})\end{equation} (cf. (1.2)).
Moreover, we have the subalgebra
\begin{equation}P(n-1)_{\bar{0}}^+=\sum_{1\leq i<j\leq
n}\mathbb{C}(E_{i,j}-E_{n+j,n+i}).\label{p756}\end{equation} Denote
\begin{equation}|\alpha|=\sum_{i=1}^n\alpha_i,\qquad x^\alpha=\prod_{i=1}^nx_i^{\alpha_i}\qquad\mbox{for}\;\alpha=
(\alpha_1,...,\alpha_n)\in\mathbb{N}^{\:n}.\end{equation}
 Set
\begin{equation}
{\cal
A}^{r}_{k,\;t}=\mbox{Span}\:\{x^\alpha\theta_{i_1}\cdots\theta_{i_t}\mid\alpha
\in\mathbb{N}^{\:n}; i_1,\cdots,i_t\in\overline{1,n};
\sum\limits_{j=r+1}^n\alpha_j-\sum\limits_{i=1}^r\alpha_i=k-t\},
\end{equation}
where $0\leq t\leq n$. Then all ${\cal A}^{r}_{k,\;t}$ form
$P(n-1)_{\bar{0}}$-submodules. Moreover,
\begin{equation}{\cal
A}^{r}_k=\bigoplus\limits_{0\leq t\leq \rm{min}\{k,n\}}{\cal
A}^{r}_{k,t}.\end{equation} Let
\begin{equation}\xi=\sum_{i=1}^nx_i\theta_i.\end{equation}
 We
first introduce a lemma, which will be used for proving linear
independence.
\begin{lem}
Let $A=(a_{i,j})$ be an $m\times m$ real matrix. If
$|a_{i,i}|>\sum\limits_{j\neq i}|a_{i,j}|$ for $i\in\ol{1,m}$, then
$det(A)\neq 0$.
\end{lem}
\psp

We study the $P(n-1)$-module structure of ${\cal A}^r_k$ according
to the following three cases.

\subsection{$r=0$}

In this case, ${\cal A}^0_k=\{0\}$ if $k<0$. Moreover, all ${\cal
A}_k^0$ are finite-dimensional $P(n-1)_{\bar 0}$-modules, which
are completely reducible. Denote by $U({\cal G})$ the universal
enveloping algebra of a Lie (super) algebra ${\cal G}$.

\begin{lem}Suppose $k>1$ and $1\leq t<min\{k,n\}$. As
$P(n-1)_{\bar{0}}$ submodules,
\begin{equation}{\cal
A}^{0}_{k,t}=U(P(n-1)_{\bar{0}})(x_1^{k-t}\theta_{n-t+1}\cdots\theta_n)\oplus
U(P(n-1)_{\bar{0}})(x_1^{k-t-1}\theta_{n-t+2}\cdots\theta_n\xi).\end{equation}
Moreover, the set
\begin{equation}
\begin{array}{l}
\{(n+k-2t)x^\alpha\theta_{i_1}\cdots\theta_{i_t}+\sum\limits_{p=1}^t\sum\limits_{s=1}^n(-1)^p\alpha_{i_p}
x^{\alpha-\epsilon_{i_p}+\epsilon_s}\theta_s\theta_{i_1}\cdots\hat{\theta_{i_p}}\cdots\theta_{i_t}\\
\mid |\alpha|=k-t;i_1,\cdots,i_t<n\mbox{ if
}\exists\alpha_{i_p}\neq0\}
\end{array}\label{p13}
\end{equation}
is a basis of
$U(P(n-1)_{\bar{0}})(x_1^{k-t}\theta_{n-t+1}\cdots\theta_n)$ and the
set
\begin{equation}
\{\sum\limits_{p=1}^t\sum\limits_{s=1}^n(-1)^p\alpha_{i_p}
x^{\alpha-\epsilon_{i_p}+\epsilon_s}\theta_s\theta_{i_1}\cdots\hat{\theta_{i_p}}\cdots\theta_{i_t}
\mid |\alpha|=k-t;\;n\in\{i_1,\cdots,i_t\}\}\label{p14}
\end{equation}
is a basis of
$U(P(n-1)_{\bar{0}})(x_1^{k-t-1}\theta_{n-t+2}\cdots\theta_n\xi)$.
\end{lem}
\begin{proof}
Denote
\begin{equation}\epsilon_i=(0,...,0,\stackrel{i}{1},0,...,0)\qquad\mbox{for}\;i\in\ol{1,n}.
\end{equation}
First we claim that
\begin{equation}
\sum\limits_{p=1}^{t+1}(-1)^p\beta_{i_p}x^{\beta-\epsilon_{i_p}}\theta_{i_1}\cdots\hat{\theta_{i_p}}\cdots\theta_{i_{t+1}}
\in U(P(n-1)_{\bar{0}})(x_1^{k-t}\theta_{n-t+1}\cdots\theta_n)
\end{equation}
for distinct $i_p\in\ol{1,n}$ and $\beta\in \mathbb{N}^{\:n}$ such
that $|\beta|=k-t+2$. It it straightforward to verify
\begin{equation}x_{i_1}^{k-t}\theta_{i_2}\cdots\theta_{i_{t+1}}\in
U(P(n-1)_{\bar{0}})(x_1^{k-t}\theta_{n-t+1}\cdots\theta_n).\end{equation}
Thus we have
\begin{equation}
x^{\alpha}\theta_{i_2}\cdots\theta_{i_{t+1}}\in
U(P(n-1)_{\bar{0}})(x_1^{k-t}\theta_{n-t+1}\cdots\theta_n)\end{equation}
for $\alpha\in \mathbb{N}^{\:n}$ such that $|\alpha|=k-t$ and
$\alpha_{i_2}=\cdots=\alpha_{i_{t+1}}=0$. Since
\begin{eqnarray}&
&\sum\limits_{p=1}^{t+1}(-1)^p\beta_{i_p}x^{\beta-\epsilon_{i_p}}\theta_{i_1}
\cdots\hat{\theta_{i_p}}\cdots\theta_{i_{t+1}}
=-\frac{\beta_{i_1}!}{(\beta_{i_1}+\cdots+\beta_{i_{t+1}}-1)!}\nonumber\\&
&\times
\prod\limits_{p=2}^{t+1}(E_{i_p,i_1}-E_{n+i_1,n+i_p})^{\beta_{i_p}}\big(x^{\beta+\sum\limits_{p=2}^{t+1}\beta_{i_p}\epsilon_{i_1}-\sum\limits_{p=2}^{t+1}\beta_{i_p}\epsilon_{i_p}-\epsilon_{i_1}}\theta_{i_2}\cdots\theta_{i_{t+1}}\big),
\end{eqnarray}
(2.13) holds.

Note
\begin{equation}A(\xi)=0\qquad\mbox{for}\;A\in P(n-1)_{\bar{0}}.\end{equation}
By the similar arguments as in the above, we have
\begin{equation}\xi(\sum\limits_{p=1}^t(-1)^p\alpha_{i_p}
x^{\alpha-\epsilon_{i_p}}\theta_{i_1}\cdots\hat{\theta_{i_p}}\cdots\theta_{i_t})\\
\in
U(P(n-1)_{\bar{0}})(x_1^{k-t-1}\theta_{n-t+2}\cdots\theta_n\xi)
\end{equation}
for $\alpha\in\mathbb{N}^{\:n}$ such that $|\alpha|=k-t$. Set
\begin{eqnarray}& &
f(\alpha;i_1,\cdots,i_t)=(k+n-2t)x^\alpha\theta_{i_1}\cdots\theta_{i_t}\nonumber\\
&
&+\sum\limits_{p=1}^t\sum\limits_{s=1}^n(-1)^p\alpha_{i_p}x^{\alpha-\epsilon_{i_p}+\epsilon_s}\theta_s\theta_{i_1}\cdots\hat{\theta_{i_p}}\cdots\theta_{i_t},
\end{eqnarray}
and
\begin{equation}
\begin{array}{lll}
g(\alpha;i_1,\cdots,i_t)&=&\xi(\sum\limits_{p=1}^t(-1)^p\alpha_{i_p}
x^{\alpha-\epsilon_{i_p}}\theta_{i_1}\cdots\hat{\theta_{i_p}}\cdots\theta_{i_t})\\
&=&\sum\limits_{p=1}^t\sum\limits_{s=1}^n(-1)^p\alpha_{i_p}x^{\alpha-\epsilon_{i_p}+\epsilon_s}\theta_s\theta_{i_1}\cdots\hat{\theta_{i_p}}\cdots\theta_{i_t}.
\end{array}\end{equation} Since
\begin{eqnarray}
f(\alpha;i_1\cdots,i_t)&=&\sum\limits_{s\neq
i_1,\cdots,i_t}\big(\sum\limits_{p=1}^t(-1)^p\alpha_{i_p}x^{\alpha+\epsilon_s-\epsilon_{i_p}}\theta_s\theta_{i_1}\cdots\hat{\theta_{i_p}}\cdots\theta_{i_t}
\nonumber \\ &
&+(\alpha_s+1)x^\alpha\theta_{i_1}\cdots\theta_{i_t}\big),
\end{eqnarray} we get $f(\alpha;i_1,\cdots,i_t)\in
P(n-1)_{\bar{0}}(x_1^{k-t}\theta_{n-t+1}\cdots\theta_n)$. Since
\begin{equation}
x^\alpha\theta_{i_1}\cdots\theta_{i_t}=\frac{1}{k+n-2t}\big(f(\alpha;i_1,\cdots,i_t)-g(\alpha;i_1,\cdots,i_t)\big),\label{p7171}
\end{equation}
we have \begin{equation}{\cal
A}^{0}_{k,t}=U(P(n-1)_{\bar{0}})(x_1^{k-t}\theta_{n-t+1}\cdots\theta_n)\oplus
U(P(n-1)_{\bar{0}})(x_1^{k-t-1}\theta_{n-t+2}\cdots\theta_n\xi).\end{equation}
Note that if $\alpha_n>0$,
\begin{equation}
f(\alpha;i_1,\cdots,i_{t-1},n)=(-1)^t\sum\limits_{s\neq
i_1,\cdots,i_{t-1},n}f(\alpha+\epsilon_s-\epsilon_n;s,i_1,\cdots,i_{t-1})
\end{equation}
and if $n\notin\{i_1,\cdots,i_t\}$,
\begin{equation}
g(\alpha;i_1,\cdots,i_t)=\frac{1}{\alpha_s+1}
\sum\limits_{p=1}^t(-1)^p\alpha_{i_p}g(\alpha+\epsilon_n-\epsilon_{i_p};n,i_1,
\cdots,\hat{i_p},\cdots,i_t).
\end{equation}
So the set (\ref{p13}) spans
$U(P(n-1)_{\bar{0}})(x_1^{k-t}\theta_{n-t+1}\cdots\theta_n)$ and
the set (\ref{p14}) spans
$U(P(n-1)_{\bar{0}})(x_1^{k-t-1}\theta_{n-t+2}\cdots\theta_n\xi)$.

We still have to check the linear independence of (\ref{p13}) and
(\ref{p14}). Assume
\begin{equation}
\sum c^\alpha_{i_1,\cdots,i_t}f(\alpha;i_1,\cdots,i_t)=0.
\end{equation}
Considering the coefficient of
$x^\alpha\theta_{i_1}\cdots\theta_{i_t}$, we obtain
\begin{equation}
(n-t+\sum\limits_{s\neq
i_1,\cdots,i_t}\alpha_s)c^\alpha_{i_1,\cdots,i_t}+\sum\limits_{\alpha_{i_p}>0,s\neq
i_1,\cdots,i_t,n}(-1)^p(\alpha_s+1)c^{\alpha-\epsilon_{i_p}+\epsilon_s}_{s,i_1,
\cdots,\hat{i_p},\cdots,i_t}=0.\end{equation}  By Lemma 2.1, we
get $c^\alpha_{i_1,\cdots,i_t}=0$. If
\begin{equation}
\sum
d^\alpha_{i_1,\cdots,i_t}g(\alpha;i_1,\cdots,i_t)=0,\end{equation}
then
\begin{equation}
-(\sum\limits_{p=1}^t\alpha_{i_p})d^\alpha_{i_1,\cdots,i_t}+\sum\limits_{\alpha_{i_p}>0,i_p\neq
n,s\neq
i_1,\cdots,i_t,n}(-1)^p(\alpha_s+1)d^{\alpha-\epsilon_{i_p}+\epsilon_s}_{s,i_1,\cdots,\hat{i_p},\cdots,i_t}=0,
\end{equation}
\begin{equation}
(-1)^p\alpha_{i_p}d^\alpha_{i_1,\cdots,i_t}-(\sum\limits_{q=1,q\neq
p}^t\alpha_{i_q}+\alpha_s+1)d^{\alpha-\epsilon_{i_p}+\epsilon_s}_{s,i_1,\cdots,\hat{i_p},\cdots,i_t}+\cdots=0,\;\;i_p\neq
n,\alpha_{i_p}>0.
\end{equation}
By Lemma 2.1 again, we get $d^{\alpha}_{i_1,\cdots,i_t}=0$.
\end{proof}
\begin{thm}
The subspace ${\cal A}^{0}_1$ is the natural irreducible
$P(n-1)$-module. Moreover, ${\cal A}^{0}_k$ ($k>1$) has a
composition series
\begin{equation}
{\cal A}^{0}_k\supset\la x_1^k\ra\supset\{0\} \mbox{ if }k\neq1,n,
\end{equation}
\begin{equation}
{\cal A}^{0}_n\supset\la \theta_1\cdots\theta_n\ra\supset\la
x_1^n\ra\supset\{0\}.
\end{equation}
In addition,
 \begin{equation}\la x_1^k\ra=\bigoplus\limits_{0\leq t\leq
k,n-1}U(P(n-1)_{\bar{0}})(x_1^{k-t}\theta_{n-t+1}\cdots\theta_n).\end{equation}
\end{thm}
\begin{proof}
The first statement is obvious. We will prove the others by
several steps.

\textsl{1) The submodule $\la x_1^k\ra\;(k>0)$ is irreducible.}

\psp Applying $E_{i,n+i}\mid_{{\cal A}^r}=x_i\partial_{\theta_i}$
to $0\neq f\in{\cal A}^{0}_k$, we have $x^\alpha\in\la f\ra$ for
some $\alpha$ with $|\alpha|=k$. Since $\mbox{Span}\{
x^\alpha\mid\alpha\in\mathbb{N}^n,\; |\alpha|=k\}$  is an
irreducible $P(n-1)_{\bar{0}}-$submodule, we get $x_1^k\in\la
f\ra$.\psp

\textsl{2) $\la x_1^k\ra=\bigoplus\limits_{0\leq t\leq n-1,k}
U(P(n-1)_{\bar{0}})(x_1^{k-t}\theta_{n-t+1}\cdots\theta_n)$ as
$P(n-1)_{\bar{0}}$-submodule.}

\psp In fact,
\begin{equation}
x_1^{k-t}\theta_{n-t+1}\cdots\theta_n=\frac{(k-t)!}{k!}(E_{2n-t+1,1}-E_{n+1,n-t+1})
\cdots(E_{2n,1}-E_{n+1,n})(x_1^k).\end{equation} Hence
\begin{equation}
U(P(n-1)_{\bar{0}})(x_1^{k-t}\theta_{n-t+1}\cdots\theta_n)\subset\la
x_1^k\ra\ \ \mbox{ for } 0\leq t\leq n-1,k.\end{equation} Now we
have to show that $\bigoplus\limits_{0\leq t\leq n-1,k}
U(P(n-1)_{\bar{0}})(x_1^{k-t}\theta_{n-t+1}\cdots\theta_n)$ is a
$P(n-1)$-submodule. It is sufficient to check
\begin{equation}
P(n-1)_{\bar{1}}(x_1^{k-t}\theta_{n-t+1}\cdots\theta_n)\subset\bigoplus\limits_{0\leq
t\leq n-1,k}
U(P(n-1)_{\bar{0}})(x_1^{k-t}\theta_{n-t+1}\cdots\theta_n)\label{p751}
\end{equation}
by PBW Theorem.

Note that if $n-t<i,j\leq n$,
\begin{eqnarray}
&&(E_{i,n+j}+E_{j,n+i})(x_1^{k-t}\theta_{n-t+1}\cdots\theta_n)\nonumber\\
=&&(-1)^{j-n+t}x_1^{k-t}x_i\theta_{n-t+1}\cdots\hat{\theta_j}
\cdots\theta_n+(-1)^{i-n+t}x_1^{k-t}x_j\theta_{n-t+1}\cdots
\hat{\theta_i}\cdots\theta_n.
\end{eqnarray}
When $1\leq i\leq n-t<j\leq n$,
\begin{equation}
(E_{i,n+j}+E_{j,n+i})(x_1^{k-t}\theta_{n-t+1}\cdots\theta_n)=
(-1)^{j-n+t}x_1^{k-t}x_i\theta_{n-t+1}\cdots\hat{\theta_j}\cdots\theta_n.
\end{equation}
In the case $1\leq i,j\leq n-t$,
\begin{equation}
(E_{i,n+j}+E_{j,n+i})(x_1^{k-t}\theta_{n-t+1}\cdots\theta_n)=0.
\end{equation}
Thus
\begin{equation}
(E_{i,n+j}+E_{j,n+i})(x_1^{k-t}\theta_{n-t+1}\cdots\theta_n)\subset
U(P(n-1)_{\bar{0}})(x_1^{k-t+1}\theta_{n-t+2}\cdots\theta_n).
\end{equation}
Since
\begin{equation}
(E_{n+i,j}-E_{n+j,i})(x_1^{k-t}\theta_{n-t+1}\cdots\theta_n)=0,\;\;\mbox{if
}i,j\neq1,
\end{equation}
and
\begin{equation}
(E_{n+i,1}-E_{n+1,i})(x_1^{k-t}\theta_{n-t+1}\cdots\theta_n)=
(k-t)x_1^{k-t-1}\theta_i\theta_{n-t+1}\cdots\theta_n,
\end{equation}
we get
\begin{equation}
(E_{n+i,j}-E_{n+j,i})(x_1^{k-t}\theta_{n-t+1}\cdots\theta_n)\subset
U(P(n-1)_{\bar{0}})(x_1^{k-t-1}\theta_{n-t}\cdots\theta_n).
\end{equation}
Hence (\ref{p751}) holds.\psp

\textsl{3) ${\cal A}^{0}_k/\la x_1^k\ra$ is irreducible when
$k\neq 1,n$.}

\psp Since
\begin{eqnarray}
{\cal A}_k^{0}=&&\bigoplus\limits_{0\leq t\leq n-1,k}
U(P(n-1)_{\bar{0}})(x_1^{k-t}\theta_{n-t+1}\cdots\theta_n)\nonumber\\
&&\bigoplus\limits_{0\leq t\leq n,k-1}
U(P(n-1)_{\bar{0}})(x_1^{k-t-1}\xi\theta_{n-t+2}\cdots\theta_n),
\end{eqnarray}
\begin{equation}
{\cal A}^{0}_k/\la x_1^k\ra\cong\bigoplus\limits_{0\leq t\leq
n,k-1}
U(P(n-1)_{\bar{0}})(x_1^{k-t-1}\xi\theta_{n-t+2}\cdots\theta_n)
\end{equation}
as $P(n-1)_{\bar{0}}$-modules. For any $ f\in{\cal
A}^{0}_k\setminus\la x_1^k\ra$, there exists
$\overline{x_1^{k-t_0-1}\xi\theta_{n-t_0+2}\cdots\theta_n}\in\la\bar{f}\ra$
for some $0\leq t_0\leq n,k-1$ by the complete reducibility of
${\cal A}_k^0$ as a $P(n-1)_{\bar{0}}$-module, where $\bar{f}$ is
the image of $f$ in ${\cal A}_k^{0}/\la x_1^k\ra$. We calculate
\begin{eqnarray}
&&\prod\limits_{j=t}^{t_0-1}(E_{2n-j+2,1}-E_{n+1,n-j+2})(\overline{x_1^{k-t_0-1}\xi\theta_{n-t_0+2}\cdots\theta_n})\nonumber\\
=&&\frac{(-1)^{t_0-t}(k-t_0-1)!}{(k-t-1)!}\overline{x_1^{k-t-1}\xi\theta_{n-t+2}\cdots\theta_n}
\end{eqnarray}
for $t<t_0$ and
\begin{eqnarray}
&&\prod\limits_{j=n-t'+2}^{n-t_0+2}(E_{1,n+j}+E_{j,n+1})(\overline{x_1^{k-t_0-1}\xi\theta_{n-t_0+2}\cdots\theta_n})\nonumber\\
=&&(-1)^{t'-t_0}\overline{x_1^{k-t'-1}\xi\theta_{n-t'+2}\cdots\theta_n}
\end{eqnarray}
for $t'>t_0$, which imply $\la \bar{f}\ra={\cal A}_k^{0}/\la
x_1^k\ra$.\psp

\textsl{4)} It is easy to check
$P(n-1)(\theta_1\cdots\theta_n)\subset\la x_1^n\ra.$ So $\la
\theta_1\cdots\theta_n\ra/\la
x_1^n\ra=\mathbb{C}\overline{\theta_1\cdots\theta_n}$ is
irreducible. By the similar argument as in 3), we get ${\cal
A}_n^{0}/\la\theta_1\cdots\theta_n\ra$ is irreducible.
\end{proof}

\subsection{${r=n}$}

In this case, ${\cal A}^{n}_k=0$ if $k>n$. Moreover, all ${\cal
A}_k^n$ are finite-dimensional $P(n-1)_{\bar 0}$-modules, which are
completely reducible. Again we deal with the
$P(n-1)_{\bar{0}}$-submodule ${\cal A}^{n}_{k,t}$ first (cf. (2.5))
\begin{lem}
When $k<t\leq n$,
\begin{eqnarray}
{\cal A}^{n}_{k,t}=&&
U(P(n-1)_{\bar{0}})(x_n^{t-k-1}\sum\limits_{p=n-t}^n(-1)^px_p\theta_{n-t}\cdots\hat{\theta_p}\cdots\theta_n)\nonumber\\
&&\oplus
U(P(n-1)_{\bar{0}})(x_n^{t-k}\theta_{n-t+1}\cdots\theta_n)\label{p6161}
\end{eqnarray}
as $P(n-1)_{\bar{0}}$-submodules. Moreover, the set
\begin{equation}
\begin{array}{l}
\{(t-k)x^\alpha\theta_{i_1}\cdots\theta_{i_t}
+\sum\limits_{p=1}^t\sum\limits_{s=1}^n(-1)^p\alpha_sx^{\alpha+\epsilon_{i_p}-\epsilon_s}\theta_s\theta_{i_1}\cdots\hat{\theta_{i_p}}\cdots\theta_{i_t}\mid
\alpha\in\mathbb{N}^{n},\\ |\alpha|=t-k;\;\exists j>i_1,\cdots,i_t
\mbox{ such that } \alpha_j\neq0\}
\end{array}\label{p17}
\end{equation}
is a basis of
$U(P(n-1)_{\bar{0}})(x_n^{t-k-1}\sum\limits_{p=n-t}^n(-1)^px_p\theta_{n-t}\cdots\hat{\theta_p}\cdots\theta_n)$
and the set
\begin{equation}
\begin{array}{l}
\{tx^\alpha\theta_{i_1}\cdots\theta_{i_t}
-\sum\limits_{p=1}^t\sum\limits_{s=1}^n(-1)^p\alpha_sx^{\alpha+\epsilon_{i_p}-\epsilon_s}\theta_s\theta_{i_1}\cdots\hat{\theta_{i_p}}\cdots\theta_{i_t}\mid
\alpha\in\mathbb{N}^{n},\\ |\alpha|=t-k;\; \alpha_j=0\mbox{ for
all }j>i_1,\cdots,i_t\}
\end{array}\label{p18}
\end{equation}
is a basis of
$U(P(n-1)_{\bar{0}})(x_n^{t-k}\theta_{n-t+1}\cdots\theta_n)$.
\end{lem}
\begin{proof}
Set
\begin{equation}
f'(\alpha;i_1,\cdots,i_t)=(t-k)x^\alpha\theta_{i_1}\cdots\theta_{i_t}
+\sum\limits_{p=1}^t\sum\limits_{s=1}^n(-1)^p\alpha_sx^{\alpha+\epsilon_{i_p}-\epsilon_s}\theta_s\theta_{i_1}\cdots\hat{\theta_{i_p}}\cdots\theta_{i_t}
\end{equation}
and
\begin{equation}
g'(\alpha;i_1,\cdots,i_t)=tx^\alpha\theta_{i_1}\cdots\theta_{i_t}
-\sum\limits_{p=1}^t\sum\limits_{s=1}^n(-1)^p\alpha_sx^{\alpha+\epsilon_{i_p}-\epsilon_s}\theta_s\theta_{i_1}\cdots\hat{\theta_{i_p}}\cdots\theta_{i_t}.
\end{equation} Denote
\begin{equation}\tilde U_t=U(P(n-1)_{\bar{0}})(x_n^{t-k-1}\sum\limits_{p=n-t}^n(-1)^px_p\theta_{n-t}
\cdots\hat{\theta_p}\cdots\theta_n).\end{equation}

It is straightforward to check
\begin{equation}x_{i_0}^{t-k-1}\sum\limits_{p=0}^t(-1)^px_{i_p}\theta_{i_0}\cdots\hat{\theta_{i_p}}\cdots\theta_t\in
\tilde U_t.\end{equation} for $i_0,\cdots,i_t\in\overline{1,n}$.
Consequently,
\begin{equation}x^\alpha\sum\limits_{p=0}^t(-1)^px_{i_p}\theta_{i_0}\cdots\hat{\theta_{i_p}}\cdots\theta_t\in
\tilde U_t.\end{equation} for $\alpha\in\mathbb{N}^{n}$ such that
$|\alpha|=t-k-1,\alpha_s=0$ for $s\neq i_0,\cdots,i_t$. Note
\begin{eqnarray}& &\prod\limits_{s\neq
i_0,\cdots,i_t}(-E_{i_0,s}+E_{n+s,n+i_0})^{\alpha_s}\big(f'(\alpha+\sum\limits_{s\neq
i_0,\cdots,i_t}\alpha_s\epsilon_{i_0}-\sum\limits_{s\neq
i_0,\cdots,i_t}\alpha_s\epsilon_{i_s};i_1,\cdots,i_t)\big)\nonumber\\
&=&\frac{(\alpha_{i_0}+\sum\limits_{s\neq
i_0,\cdots,i_t}\alpha_s-1)!}{\alpha_{i_0}!}f'(\alpha;i_1,\cdots,\i_t).
\end{eqnarray}
So $f'(\alpha;i_1,\cdots,i_t)\in \tilde U_t$ for all
$\alpha\in\mathbb{N}^{n}$ with $|\alpha|=t-k$.

Since $x_{i_1}^{t-k}\theta_{i_1}\cdots\theta_{i_t}\in
U(P(n-1)_{\bar{0}})(x_n^{t-k}\theta_{n-t+1}\cdots\theta_n),$ we
have $x^\alpha\theta_{i_1}\cdots\theta_{i_t}\in
U(P(n-1)_{\bar{0}})(x_n^{t-k}\theta_{n-t+1}\cdots\theta_n)$ for
any $\alpha\in\mathbb{N}^{n}$ with $|\alpha|=t-k$ and $\alpha_s=0$
if $s\neq i_1,\cdots,i_t$. Note
\begin{eqnarray}
g'(\alpha;i_1,\cdots,i_t)&=&\sum\limits_{p=1}^t\frac{(\alpha_{i_p}+1)!}{(\alpha_{i_p}+\sum\limits_{s\neq
i_1,\cdots,i_t}\alpha_s)!}\prod\limits_{s\neq
i_1,\cdots,i_t}(-E_{i_p,s}+E_{n+s,n+i_p})^{\alpha_s}\nonumber
\\&&(x^{\alpha+\sum\limits_{s\neq
i_1,\cdots,i_t}\alpha_s\epsilon_{i_p}-\sum\limits_{s\neq
i_1,\cdots,i_t}\alpha_s\epsilon_s}\theta_{i_1}\cdots\theta_{i_t}).
\end{eqnarray}
So we obtain $g'(\alpha;i_1,\cdots,i_t)\in
U(P(n-1)_{\bar{0}})(x_n^{t-k}\theta_{n-t+1}\cdots\theta_n)$ for
all $\alpha\in\mathbb{N}^{n}$ with $|\alpha|=t-k$. On the other
hand,
\begin{eqnarray}
{\cal
A}_{k,t}^{n}&=&\mbox{Span}\;\{f'(\alpha;i_1,\cdots,i_t),g'(\alpha,i_t,\cdots,i_t)
\mid\alpha\in\mathbb{N}^{n},\nonumber
\\ & &|\alpha|=t-k;\;i_1,\cdots,i_t\in\overline{1,n}\},
\end{eqnarray}
which implies (\ref{p6161}). If $\alpha_s=0$ for
$s>i_1,\cdots,i_t$, we can assume $i_t=\mbox{max}\{j\mid
\alpha_j>0\}$ and have
\begin{equation}
f'(\alpha;i_1,\cdots,i_t)=\frac{(-1)^t}{\alpha_{i_t}+1}\sum\limits_{s\in\overline{1,n}\setminus\{i_1,\cdots,i_t\}}\alpha_sf'(\alpha+\epsilon_{i_t}-\epsilon_s;s,i_1,\cdots,i_{t-1}).
\label{p20}
\end{equation}
If $\exists j>i_1,\cdots,i_t$ such that $\alpha_j>0$, let
$s=\mbox{max}\{j\mid \alpha_j>0\}$. Note  that $s>i_1,\cdots,i_t$
and
\begin{equation}
g'(\alpha;i_1,\cdots,i_t)=-\sum\limits_{p=1}^t(-1)^pg'(\alpha+\epsilon_{i_p}-\epsilon_s;s,i_1,\cdots,\hat{i_p},\cdots,i_t).\label{p21}
\end{equation}
Now (2.48) spans $\tilde U_t$ and (2.49) spans
$U(P(n-1)_{\bar{0}})(x_n^{t-k}\theta_{n-t+1}\cdots\theta_n)$ by
(3.58)-(3.59). Lemma 2.1 implies that the sets (\ref{p17}) and
(\ref{p18}) are linear independent.
\end{proof}
\begin{thm} The subspace
${\cal A}^{n}_k=0$ when $k>n$. Moreover, the submodules ${\cal
A}^{n}_n$ and ${\cal A}^{n}_{n-1}$ are irreducible. When $k\leq
n-2$, we have the following composition series:
\begin{equation}
{\cal A}^{n}_k\supset\la
x_n^{n-k-2}\sum\limits_{i=1}^n(-1)^ix_i\theta_1\cdots\hat{\theta_i}\cdots\theta_n\ra
\supset\{0\}\qquad\mbox{ if }k\leq n-2,k\neq0,
\end{equation}
\begin{equation}
{\cal A}^{n}_0\supset\la1\ra\supset\la
x_n^{n-2}\sum\limits_{i=1}^n(-1)^ix_i\theta_1\cdots\hat{\theta_i}\cdots\theta_n\ra\supset\{0\}.
\end{equation}
Furthermore, in terms of (2.52),
\begin{equation}
\la
x_n^{n-k-2}\sum\limits_{i=1}^n(-1)^ix_i\theta_1\cdots\hat{\theta_i}\cdots\theta_n\ra
=\bigoplus_{max\{0,k+1\}\leq t\leq n}\tilde U_t.\label{p6162}
\end{equation}
\end{thm}
\begin{proof}
The first and second statement are trivial. Now assume $k\leq
n-2$.

\textsl{1) The module $\la
x_n^{n-k-2}\sum\limits_{i=1}^n(-1)^ix_i\theta_1\cdots\hat{\theta_i}\cdots\theta_n
\ra$ is the minimal submodule of ${\cal A}^{n}_k$.}

Given $0\neq f\in{\cal A}_k^{n}$. There should be a weight vector
$g\in\la f\ra$ such that $P(n-1)_{\bar{0}}^+(g)=0$ by the completely
reducibility of ${\cal A}_k^n$ as a $P(n-1)_{\bar 0}$-module. Up to
a scalar multiple,
\begin{equation}g=x_n^{n-s-1-k}\sum\limits_{p=s}^n(-1)^px_p\theta_s\cdots
\hat{\theta_p}\cdots\theta_n\;\mbox{or}\;
x_n^{n-s+1-k}\theta_s\cdots\theta_n\end{equation} for some $s$. If
$g=x_n^{n-s-1-k}\sum\limits_{p=s}^n(-1)^px_p\theta_s\cdots\hat{\theta_p}
\cdots\theta_n$,
then
\begin{equation}
\prod\limits_{i=1}^{s-1}(E_{i,2n}-E_{n,n+i})(g)
=x_n^{n-k-2}\sum\limits_{i=1}^n(-1)^ix_i\theta_1\cdots\hat{\theta_i}\cdots\theta_n\in\la
f\ra.
\end{equation}
When $g=x_n^{n-s+1-k}\theta_s\cdots\theta_n$, we have
\begin{eqnarray}
x_n^{n-k}\theta_1\cdots\theta_n&=&\prod\limits_{i=1}^{s-1}(E_{2n,i}-E_{n+i,n})(g)\in\la
f\ra,\\
x_n^{n-k-1}\theta_1\cdots\theta_{n-1}&=&\frac{(-1)^{n-1}}{n-k}E_{n,2n}.(x_n^{n-k}\theta_1\cdots\theta_n)\in\la
f\ra,\\
x_px_n^{t-k-1}\theta_1\cdots\hat{\theta_p}\cdots\theta_n&=&\frac{(-1)^{p-1}}{n-k-1}\big((E_{p,2n}+E_{n,n+p})(x_px_n^{n-k-1}\theta_1\cdots\theta_n)\nonumber\\
&&+(-1)^nx_n^{n-k-1}\theta_1\cdots\theta_{n-1}\big).
\end{eqnarray}
Anyway, $
x_n^{n-k-2}\sum\limits_{i=1}^n(-1)^ix_i\theta_1\cdots\hat{\theta_i}\cdots\theta_n\in\la
f\ra$.\psp

\textsl{2) Equation (\ref{p6162}) holds.}

\psp Since
\begin{eqnarray}
&&x_n^{t-k-1}\sum\limits_{p=n-t}^n(-1)^px_p\theta_{n-t}\cdots\hat{\theta_p}\cdots\theta_n\nonumber\\
=&&\frac{(t-k-1)!}{(n-k-2)!}\prod\limits_{i=1}^{n-t-1}(-E_{i,2n}-E_{n,n+i})(x_n^{n-k-2}\sum\limits_{i=1}^n(-1)^ix_i\theta_1\cdots\hat{\theta_i}\cdots\theta_n),
\end{eqnarray}
we get
\begin{equation}
\tilde U_t\subset \la
x_n^{n-k-2}\sum\limits_{i=1}^n(-1)^ix_i\theta_1\cdots\hat{\theta_i}\cdots\theta_n\ra.
\end{equation}
It is straightforward to verify
\begin{equation}
(E_{i,n+j}+E_{j,n+i})(x_n^{t-k-1}\sum\limits_{p=n-t}^n(-1)^px_p\theta_{n-t}\cdots
\hat{\theta_p}\cdots\theta_n)\in\tilde U_{t-1}
\end{equation}
and
\begin{equation}
(E_{n+i,j}-E_{n+j,i})(x_n^{t-k-1}\sum\limits_{p=n-t}^n(-1)^px_p\theta_{n-t}\cdots\hat{\theta_p}\cdots\theta_n)\in
\tilde U_{t+1}.\end{equation} Therefore, the right side of
(\ref{p6162}) is a $P(n-1)$-submodule.\psp

\textsl{3) The quotient module ${\cal A}^{n}_k/\la
x_n^{n-k-2}\sum\limits_{i=1}^n(-1)^ix_i\theta_1\cdots\hat{\theta_i}\cdots\theta_n\ra$
is irreducible when $k\neq0$.}

\psp Since
\begin{equation}
\tilde{U_t}\subset\la
x_n^{n-k-2}\sum\limits_{i=1}^n(-1)^ix_i\theta_1\cdots\hat{\theta_i}\cdots\theta_n\ra,
\end{equation}
we have
\begin{eqnarray}& &
{\cal A}^{n}_k/\la
x_n^{n-k-2}\sum\limits_{i=1}^n(-1)^ix_i\theta_1\cdots\hat{\theta_i}\cdots\theta_n\ra
\nonumber \\& \cong&\bigoplus_{1\leq s\leq
n,n+1-k}U(P(n-1)_{\bar{0}})(x_n^{n-s+1-k}\theta_s\cdots\theta_n)
\end{eqnarray}
as $P(n-1)_{\bar{0}}$-modules. We use $\bar u$ to denote the image
in the quotient space for  $u\in{\cal A}^{n}_k$. For any $0\neq
f\in{\cal A}^{n}_k\setminus\la
x_n^{n-k-2}\sum\limits_{i=1}^n(-1)^ix_i\theta_1\cdots\hat{\theta_i}\cdots\theta_n\ra$,
we have
$\bar{g}=\overline{x_n^{n-s+1-k}\theta_s\cdots\theta_n}\in\la\bar{f}\ra$
for some $s$. Since
\begin{equation}
(E_{n,n+s}+E_{s,2n})(\bar{g})=(n-s+1-k)\overline{x_n^{n-s-k}\theta_{s+1}\cdots
\theta_n}\end{equation} and
\begin{equation}
(E_{2n,s-1}-E_{n+s-1,n})(\bar{g})=\overline{x_n^{n-s+2-k}\theta_{s-1}\cdots\theta_n},
\end{equation}
we have $\overline{x_n^{n-s+1-k}\theta_s\cdots\theta_n}\in\la
\bar{f}\ra$ for all $0<s\leq \mbox{min}\{n,n+1-k\}$. Thus $\la
\bar{f}\ra={\cal A}^{n}_k/\la
x_n^{n-k-2}\sum\limits_{i=1}^n(-1)^ix_i\theta_1\cdots\hat{\theta_i}\cdots\theta_n\ra$.\psp

\textsl{4) $k=0$.}\psp

Since
\begin{eqnarray}
&&(E_{i,j}-E_{n+j,n+i})(1)=(-x_j\partial_{x_i}-\theta_j\partial_{\theta_i})(1)=0,\\
&&(E_{i,n+j}+E_{j,n+i})(1)=(\partial_{x_i}\partial_{\theta_j}+\partial_{x_j}
\partial_{\theta_i})(1)=0\\
and &&(E_{n+i,j}-E_{n+j,i})(1)=x_i\theta_j-x_j\theta_i\in\la
x_n^{n-k-2}\sum\limits_{i=1}^n(-1)^ix_i\theta_1\cdots\hat{\theta_i}\cdots\theta_n
\ra
\end{eqnarray}
for $i,j\in\ol{1,n}$, we get
\begin{equation}
P(n-1)(1)\subset\la
x_n^{n-k-2}\sum\limits_{i=1}^n(-1)^ix_i\theta_1\cdots\hat{\theta_i}\cdots\theta_n
\ra
\end{equation}
by (2.48) with $t=1$, (2.52) and (2.62). Thus $\la 1\ra/\la
x_n^{n-k-2}\sum\limits_{i=1}^n(-1)^ix_i\theta_1\cdots\hat{\theta_i}\cdots\theta_n
\ra=\mathbb{C}\bar{1}$ is irreducible. By the similar arguments as
those in 3), we can prove that ${\cal A}^{n}_0/\la 1\ra$ is
irreducible.
\end{proof}

\subsection{$0<r<n$}

Set
\begin{eqnarray}
L_1&=&\sum\limits_{i=1}^{r-1}\mathbb{C}(E_{i,i}-E_{n+i,n+i}-E_{i+1,i+1}+E_{n+i+1,n+i+1})\nonumber\\
&&+\sum\limits_{i,j\in\overline{1,r},i\neq j}\mathbb{C}(E_{i,j}-E_{n+j,n+i}),\\
L_2&=&\sum\limits_{i=r+1}^{n-1}\mathbb{C}(E_{i,i}-E_{n+i,n+i}-E_{i+1,i+1}+E_{n+i+1,n+i+1})\nonumber\\
&&+\sum\limits_{i,j\in\overline{r+1,n},i\neq
j}\mathbb{C}(E_{i,j}-E_{n+j,n+i}),
\end{eqnarray}
and
\begin{equation}
L_1^+=\sum\limits_{1\leq i<j\leq
r}\mathbb{C}(E_{i,j}-E_{n+j,n+i}),\ \ L_2^+=\sum\limits_{r+1\leq
i<j\leq n}\mathbb{C}(E_{i,j}-E_{n+j,n+i}).
\end{equation}
Recall that
\begin{equation*}
{\cal A}^{r}_{k,l}=\mbox{Span}\{
x^\alpha\theta_{i_1}\cdots\theta_{i_l}\mid
i_1,\cdots,i_l\in\overline{1,n};\;\sum\limits_{j=r+1}^n\alpha_j-\sum\limits_{i=1}^r\alpha_i+l=k\}
\end{equation*}
is a $P(n-1)_{\bar{0}}$ submodule. Let $V_{k,l}^{r}$ be the
$P(n-1)_{\bar{0}}$-submodule of ${\cal A}_{k,l}^{r}$ generated by
\begin{center}{\bf Table 1}\end{center}
\begin{center}
\begin{tabular}{|c|c|}
  \hline
$x_{r+1}^{k-l}\theta_{n-l+1}\cdots\theta_n$ & if $l<n-r,l\leq
k$\\
\hline
 $x_{r}^{l-k}\theta_{n-l+1}\cdots\theta_n$ & if $k<l<n-r$\\
\hline
$\begin{array}{l}(k-l+1)x_{r+1}^{k-l}\theta_{n-l}\cdots\theta_r\theta_{r+2}\cdots\theta_n
+\\
\sum\limits_{p=n-l}^r(-1)^{p-r}x_p\theta_{n-l}\cdots\hat{\theta_p}\cdots\theta_rx_{r+1}^{k-l+1}\theta_{r+1}\cdots\theta_n\end{array}$
& if $n-r\leq l\leq k,l\neq\frac{1}{2}(k+n-r)$\\
\hline
$\sum\limits_{p=n-l-1}^r(-1)^px_p\theta_{n-l-1}\cdots\hat{\theta_p}\cdots\theta_rx_{r+1}^{k-l+1}\theta_{r+2}\cdots\theta_n$
& if $n-r\leq l\leq k,l=\frac{1}{2}(k+n-r)$\\
\hline
$x_r^{l-k-1}\sum\limits_{p=n-l}^r(-1)^{p}x_p\theta_{n-l}\cdots\hat{\theta_p}\cdots\theta_r\theta_{r+1}\cdots\theta_n$
& if $l>k,l\geq n-r$\\
\hline
\end{tabular}
\end{center}
and let $W^{r}_{k,l}$ be $P(n-1)_{\bar{0}}$-submodule generated by
\begin{center}{\bf Table 2}\end{center}
\begin{center}
\begin{tabular}{|c|c|}
\hline $x_{r+1}^{k-l-1}\xi_r\theta_{n-l+2}\cdots\theta_n$ & if
$l<k,l\leq n-r$\\
\hline

$(l-k+1)x_r^{l-k}\theta_r\theta_{n-l+2}\cdots\theta_n+x_r^{l-k+1}\xi_r\theta_{n-l+2}\cdots\theta_n$
& if $k\leq l\leq n-r,l\neq\frac{1}{2}(k+n-r)$\\
\hline

$x_r^{l-k+1}\theta_r\xi_r\theta_{n-l+3}\cdots\theta_n$ & if $k\leq
l\leq n-r,l=\frac{1}{2}(k+n-r)$\\
\hline

$x_{r+1}^{l-k}\theta_{n-l+1}\cdots\theta_n$ & if $n-r\leq l<
k$\\
\hline

$x_r^{l-k}\theta_{n-l+1}\cdots\theta_n$ & if $l\geq k,l>n-r$\\
\hline
\end{tabular}
\end{center}
where $\xi_r=\sum\limits_{i=r+1}^rx_i\theta_i$.

\begin{lem} The subspace
${\cal A}^{r}_{k,l}=V^{r}_{k,l}\oplus W^{r}_{k,l}$ if
$l\neq\frac{1}{2}(k+n-r)$. When $l=\frac{1}{2}(k+n-r)$, we have
the following composition series of $P(n-1)_{\bar 0}$-submodules:
\begin{equation}
{\cal A}^{r}_{k,l}\supset V^{r}_{k,l}+W^{r}_{k,l}\supset
V^{r}_{k,l}\;(or \ \ W^{r}_{k,l})\supset V^{r}_{k,l}\bigcap
W^{r}_{k,l}\supset\{0\}.\label{p6152}
\end{equation}
\end{lem}
\begin{proof}
Set
\begin{eqnarray}
{\cal
A}^{r}_{k',k'',t,s}&=&\mbox{Span}\{x^\alpha\theta_{i_1}\cdots\theta_{i_t}\theta_{j_1}\cdots\theta_{j_s}\mid
i_1,\cdots,i_t\in\overline{1,r},j_1,\cdots,j_s\in\overline{r+1,n};\nonumber\\&
& \alpha\in\mathbb{N}^n,\;\sum\limits_{i=1}^r\alpha_i=t-k',
\sum\limits_{i=r+1}^n\alpha_i=k''-s\}.\label{p7161}
\end{eqnarray}
Note that ${\cal
A}^{r}_{k,l}=\bigoplus\limits_{k'+k''=k,s+t=l}{\cal
A}^{r}_{k',k'',t,s}$. Moreover,
\begin{eqnarray}& &
{\cal
A}^{r}_{k',k'',t,s}=U(L_1+L_2)(x_r^{t-k'}\theta_{r-t+1}\cdots\theta_r
x_{r+1}^{k''-s-1}\xi_r\theta_{n-s+2}\cdots\theta_n)\nonumber\\
&&\oplus
U(L_1+L_2)(x_r^{t-k'}\theta_{r-t+1}\cdots\theta_rx_{r+1}^{k''-s}\theta_{n-s+1}
\cdots\theta_n)\nonumber\\
&&\oplus
U(L_1+L_2)(x_r^{t-k'-1}\sum\limits_{p=r-t}^r(-1)^px_p\theta_{r-t}\cdots
\hat{\theta_p}\cdots\theta_rx_{r+1}^{k''-s-1}\xi_r\theta_{n-s+2}\cdots\theta_n)
\nonumber\\
&&\oplus
U(L_1+L_2)(x_r^{t-k'-1}\sum\limits_{p=r-t}^r(-1)^px_p\theta_{r-t}\cdots\hat{\theta_p}\cdots\theta_rx_{r+1}^{k''-s}\theta_{n-s+1}\cdots\theta_n)\label{p7162}
\end{eqnarray}
if $0\leq s<n-r$, $0\leq t<r$. We claim that
\begin{equation}
x_r^{\alpha_r-1}\sum\limits_{p=r-t}^r(-1)^px_p\theta_{r-t}\cdots\hat{\theta_p}\cdots\theta_rx_{r+1}^{\alpha_{r+1}}\theta_{n-s+1}\cdots\theta_n\in
V^{r}_{k,l}\label{p6151}
\end{equation}
for all $0<l\leq n$, $0\leq t<r$, $0<s<n-r$, $s+t=l$ and
$\alpha_{r+1}-\alpha_r+l=k$. We prove it case by case.\psp

{\it (a) $l<n-r$}.\psp

We have
\begin{equation}
x_r^{\alpha_r}x_{r+1}^{\alpha_{r+1}}\theta_{n-l+1}\cdots\theta_n\in
V^{r}_{k,l}
\end{equation}
because
\begin{equation}
(E_{r+1,r}-E_{n+r,n+r+1})|_{{\cal
A}^{r}}=-x_rx_{r+1}+\theta_r\partial_{\theta_{r+1}}.
\end{equation}
Assume
\begin{equation}
x_r^{\alpha_r-1}\sum\limits_{p=r-t+1}^r(-1)^px_p\theta_{r-t+1}\cdots\hat{\theta_p}\cdots\theta_rx_{r+1}^{\alpha_{r+1}}\theta_{n-s}\cdots\theta_n\in
V^{r}_{k,l}.
\end{equation}
We have
\begin{equation}
h_i=x_r^{\alpha_r-1}\sum\limits_{i\neq
p=r-t}^r(-1)^p\mbox{sgn}(p-i)x_p\theta_{r-t}\cdots\hat{\theta_i}\cdots\hat{\theta_p}\cdots\theta_rx_{r+1}^{\alpha_{r+1}}\theta_{n-s}\cdots\theta_n\in
V^{r}_{k,l},
\end{equation}
where
\begin{equation}
\mbox{sgn}(p-i)=\left\{
\begin{array}{lll}
1 &\rm{if} &p>i,\\
0 &\rm{if} &p=i,\\
-1 &\rm{if} &p<i.
\end{array}
\right.
\end{equation}
Note
\begin{eqnarray}
&&\sum\limits(-1)^i(-E_{n-s,i}+E_{n+i,2n-s})(h_i)\nonumber\\
=&&-tx_r^{\alpha_r-1}\sum\limits_{p=r-t}^r(-1)^{p-r}x_p\theta_{r-t}\cdots
\hat{\theta_p}\cdots\theta_rx_{r+1}^{\alpha_{r+1}}\theta_{n-s+1}\cdots\theta_n.
\end{eqnarray}
So (\ref{p6151}) holds when $l<n-r$.\psp

{\it (b) $n-r\leq l\leq k$.}\psp

If $l\neq\frac{1}{2}(k+n-r)$, we set
\begin{eqnarray}
g_{n-l-1}&=&(-E_{r+1,n-l-1}+E_{n-l-1,r+1}).\big((k-l+1)x_{r+1}^{k-l}\theta_{n-l}\cdots\theta_r\theta_{r+2}\cdots\theta_n\nonumber\\
&&+
\sum\limits_{p=n-l}^r(-1)^{p-r}x_p\theta_{n-l}\cdots\hat{\theta_p}\cdots\theta_rx_{r+1}^{k-l+1}\theta_{r+1}\cdots\theta_n\big)\nonumber\\
&=&(k-l+1)x_{n-l-1}\theta_{n-l}\cdots\theta_rx_{r+1}^{k-l+1}\theta_{r+2}\cdots\theta_n\nonumber\\
&&+\sum\limits_{p=n-l}^r(-1)^{p-r}x_{n-l-1}x_p\theta_{n-l}\cdots\hat{\theta_p}\cdots\theta_rx_{r+1}^{k-l+2}\theta_{r+1}\cdots\theta_n\nonumber\\
&&+\sum\limits_{p=n-l}^r(-1)^{p-n+l}x_p\theta_{n-l-1}\cdots\hat{\theta_p}\cdots\theta_rx_{r+1}^{k-l+2}\theta_{r+2}\cdots\theta_n,\\
g_i&=&(-1)^{i-n+l}(-E_{r+1,i}+E_{n+i,n+r+1})(E_{i,n-l-1}E_{2n-l-1,n+i})\big((k-l+1)\theta_{n-l}\cdots\theta_r\nonumber\\
&&\times x_{r+1}^{k-l}\theta_{r+2}\cdots\theta_n +
\sum\limits_{p=n-l}^r(-1)^{p-r}x_p\theta_{n-l}\cdots\hat{\theta_p}\cdots\theta_rx_{r+1}^{k-l+1}\theta_{r+1}\cdots\theta_n\big)\nonumber\\
&=&(k-l+1)x_i\theta_{n-l-1}\cdots\hat{\theta_i}\theta_rx_{r+1}^{k-l+1}\theta_{r+2}\cdots\theta_n\nonumber\\
&&+\sum\limits_{i\neq p=n-l-1}^r(-1)^{p-r}\mbox{sgn}(p-i)x_ix_p\theta_{n-l-1}\cdots\hat{\theta_i}\cdots\hat{\theta_p}\cdots\theta_rx_{r+1}^{k-l+2}\theta_{r+1}\cdots\theta_n\nonumber\\
&&+\sum\limits_{i\neq
p=n-l-1}^r(-1)^{p+i+1}x_p\theta_{n-l-1}\cdots\hat{\theta_p}\cdots\theta_rx_{r+1}^{k-l+2}
\theta_{r+2}\cdots\theta_n
\end{eqnarray}
for $i=n-l,\cdots,r$. We calculate
\begin{eqnarray}& &
\sum\limits_{p=n-l-1}^r(-1)^px_p\theta_{n-l-1}\cdots\hat{\theta_p}
\cdots\theta_rx_{r+1}^{k-l+1}\theta_{r+2}\cdots\theta_n\nonumber\\
&=&\frac{1}{k+n-r-2l}\sum_{i=n-l-1}^r(-1)^ig_i\in V^{r}_{k,l}.
\end{eqnarray}
Again by induction on $t$, we obtain (\ref{p6151}) holds for
$n-r\leq l\leq k$. It can be similarly proved when $l>k$ and
$l\geq n-r$.

Now let
\begin{eqnarray}
f_j&=&\frac{1}{\alpha_{r+1}+1}(E_{r-t,j}-E_{n+j,n+r-t})(E_{j,r+1}-E_{n+j,n+r-t})(-E_{r+1,r}-E_{n+r,n+r+1}).\nonumber\\
&&\big(x_r^{\alpha_r-1}
\sum\limits_{p=r-t}^r(-1)^px_p\theta_{r-t}\cdots\hat{\theta_p}\cdots\theta_rx_{r+1}^{\alpha_{r+1}}\theta_{n-s+1}\cdots\theta_n\big)\nonumber\\
&=&(-1)^{r-t}x_r^{\alpha_r}\theta_{r-t+1}\cdots\theta_rx_{r+1}^{\alpha_{r+1}}\theta_{n-s+1}\cdots\theta_n\nonumber\\
&&+x_r^{\alpha_r}\sum\limits_{p=r-t+1}^r(-1)^{p-t}x_p\theta_{r-t+1}\cdots\hat{\theta_p}\cdots\theta_rx_{r+1}^{\alpha_{r+1}}x_j\theta_j\theta_{n-s+1}\cdots\theta_n,
\end{eqnarray}
$j=r+1,\cdots,n-s$. Taking summation on $j$, we get
\begin{equation}
\begin{array}{l}
(\alpha_{r+1}+n-r-s)x_r^{\alpha_r}\theta_{r-t+1}\cdots\theta_rx_{r+1}^{\alpha_{r+1}}\theta_{n-s+1}\cdots\theta_n
+x_r^{\alpha_r}\sum\limits_{p=r-t+1}^r(-1)^{p-r}x_p\\
\times\theta_{r-t+1}\cdots\hat{\theta_p}\cdots\theta_r
 x_{r+1}^{\alpha_{r+1}}\xi_r\theta_{n-s+1}\cdots\theta_n\in
V^{r}_{k,l}.\label{p7165}
\end{array}
\end{equation}
It is not difficult to verify that the subspace
\begin{eqnarray}& &
\bigoplus\limits_{\stackrel{s+t=l,}{\alpha_{r+1}-\alpha_r=k-l}}
\big\{U(L_1+L_2)(x_r^{\alpha_r-1}\sum\limits_{p=r-t}^p(-1)^px_p\theta_{r-t}\cdots
\hat{\theta_p}\cdots\theta_rx_{r+1}^{\alpha_{r+1}}\theta_{n-s+1}\cdots\theta_n)
\nonumber\\ & &\oplus U(L_1+L_2)
[(\alpha_{r+1}+n-r-s)x_r^{\alpha_r}\theta_{r-t+1}\cdots\theta_rx_{r+1}^{\alpha_{r+1}}\theta_{n-s+1}\cdots\theta_n\nonumber\\
&&+x_r^{\alpha_r}\sum\limits_{p=r-t+1}^r(-1)^{p-r}x_p\theta_{r-t+1}\cdots\hat{\theta_p}\cdots\theta_r
 x_{r+1}^{\alpha_{r+1}}\xi_r\theta_{n-s+1}\cdots\theta_n]\big\}
\end{eqnarray}
is invariant under $P(n-1)_{\bar{0}}$, which implies that it is
equal to $V^{r}_{k,l}$ exactly.

By a similar argument, we obtain that
\begin{eqnarray}
W^{r}_{k,l}&=&\bigoplus\limits_{\stackrel{s+t=l,}{\alpha_{r+1}-\alpha_r=k-l}}
\big\{U(L_1+L_2)(x_r^{\alpha_r}\theta_{r-t+1}\cdots\theta_rx_{r+1}^{\alpha_{r+1}-1}
\xi_r\theta_{n-s+2}\cdots\theta_n)\nonumber\\
&&\oplus
U(L_1+L_2)[(\alpha_r+t)x_r^{\alpha_r}\theta_{r-t+1}\cdots\theta_rx_{r+1}^{\alpha_{r+1}}\theta_{n-s+1}\cdots\theta_n+\nonumber\\
&&x_r^{\alpha_r}\sum\limits_{p=r-t+1}^r(-1)^{p-r}x_p\theta_{r-t+1}\cdots\hat{\theta_p}\cdots\theta_r
x_{r+1}^{\alpha_{r+1}}\xi_r\theta_{n-s+1}\cdots\theta_n]\big\}.
\end{eqnarray}
Now it is easy to see ${\cal A}^{r}_{k,l}=V^{r}_{k,l}\oplus
W^{r}_{k,l}$ if $l\neq \frac{1}{2}(k+n-r)$.\psp

Now assume $l=\frac{1}{2}(k+n-r)$.\psp

According to the above arguments,
\begin{eqnarray}&&
{\cal A}^{r}_{k,l}/(V^{r}_{k,l}+W^{r}_{k,l})
=\bigoplus\limits_{\stackrel{s+t=l,}{\alpha_{r+1}-\alpha_r+l=k}}U(L_1+L_2)(x_r^{\alpha_r}\theta_{r-t+1}\cdots\theta_rx_{r+1}^{\alpha_{r+1}}\theta_{n-s+1}\cdots\theta_n)\nonumber \\
&=&\bigoplus\limits_{\stackrel{s+t=l,}{\alpha_{r+1}-\alpha_r+l=k}}U(L_1+L_2)
(x_r^{\alpha_r}\sum\limits_{p=r-t+1}^r(-1)^{p}x_p\theta_{r-t+1}\cdots\hat{\theta_p}
\cdots\theta_r\nonumber\\
& &\times
x_{r+1}^{\alpha_{r+1}}\xi_r\theta_{n-s+1}\cdots\theta_n).
\end{eqnarray}

For any $0\neq f\in {\cal
A}^{r}_{k,l}\setminus(V^{r}_{k,l}+W^{r}_{k,l})$, there should be
some weight vector $\bar{g}\in U(P(n-1)_{\bar{0}})(\bar{f})$ such
that $L_1^+(\bar{g})=0$ and $L_2^+(\bar{g})=0$. Up to a scalar
multiple,
\begin{equation} g\equiv
x_r^{k_r}\sum\limits_{p=r-t_0+1}^r(-1)^{p}x_p\theta_{r-t_0+1}\cdots\hat{\theta_p}
\cdots\theta_rx_{r+1}^{k_{r+1}}\xi_r\theta_{n-s_0+1}\cdots\theta_n\;
(\mbox{mod}\;V^{r}_{k,l}+W^{r}_{k,l})
\end{equation}
for some $s_0,t_0,k_r,k_{r+1}\in\mathbb{N}$ such that $s_0+t_0=l$,
$k_{r+1}-k_r=k-l$ and $0\leq s_0<n-r$.\psp

{\it (1) $l>n-r$.}\psp

 We have
\begin{eqnarray}&&(E_{n-l-1,r+2}-E_{n+r+2,2n-l-1})\cdots(E_{r-t_0+1,n-s_0}-E_{2n-s_0,n+r-t_0+1})(g)\nonumber\\
&\equiv&(-1)^{\sum\limits_{j=r-n+l+2}^{t_0}j}\big(x_r^{k_r}
\sum\limits_{p=n-l}^r(-1)^{p}x_p\theta_{n-l}\cdots\hat{\theta_p}\cdots\theta_rx_{r+1}^{k_{r+1}+1}\theta_{r+1}
\cdots\theta_n\big)\nonumber\\ & &(\mbox{mod}
\;V^{r}_{k,l}+W^{r}_{k,l}).
\end{eqnarray}
Furthermore,
\begin{equation}
x_r^{\alpha_r}\sum\limits_{p=n-l}^r(-1)^{p}x_p\theta_{n-l}\cdots\hat{\theta_p}\cdots\theta_rx_{r+1}^{\alpha_{r+1}+1}\theta_{r+1}\cdots\theta_n
\in U(P(n-1)_{\bar{0}})(f)+V^{r}_{k,l}+W^{r}_{k,l}
\end{equation}
for all $\alpha\in\mathbb{N}^n$ such that
$\alpha_{r+1}-\alpha_r=k-l$ because
\begin{eqnarray}
(E_{r,r+1}-E_{n+r+1,n+r})|_{{\cal
A}^{r}}=\partial_{x_r}\partial_{x_{r+1}}-\theta_{r+1}\partial_{\theta_r},\\
(E_{r+1,r}-E_{n+r,n+r+1})|_{{\cal
A}^{r}}=-x_rx_{r+1}-\theta_r\partial_{\theta_{r+1}}.
\end{eqnarray}
Therefore,
\begin{equation}
{\cal A}^{r}_{k',k'',l-n+r,n-r}\subset
U(P(n-1)_{\bar{0}})(f)+V^{r}_{k,l}+W^{r}_{k,l}\qquad\mbox{ for
}k''-k'+l=k.\end{equation} By induction on $t$, we obtain
\begin{equation}
{\cal A}^{r}_{k',k'',t,s}\subset
U(P(n-1)_{\bar{0}})(f)+V^{r}_{k,l}+W^{r}_{k,l}\;\mbox{ for
}k''-k'+l=k,s+t=l,t\geq l-n+r.
\end{equation}

{\it (2) $l\leq n-r$.}\psp

 Note
\begin{eqnarray}&
&(E_{r-1,n-l+2}-E_{2n-l+2,n+r-1})\cdots(E_{r-t_0+1,n-s_0}-E_{2n-s_0,n+r-t_0+1})(g)\nonumber
\\ &\equiv&(-1)^{\sum\limits_{j=1}^{t_0}j}\big(x_r^{k_r}x_{r+1}^{k_{r+1}}\xi_r\theta_{n-l+2}\cdots\theta_n\big)\;(\mbox{mod}
\;V^{r}_{k,l}+W^{r}_{k,l}).
\end{eqnarray}
Consequently,
\begin{equation}
x_r^{\alpha_r}x_{r+1}^{\alpha_{r+1}}\xi_r\theta_{n-l+2}\cdots\theta_n\in
U(P(n-1)_{\bar{0}})(f)+V^{r}_{k,l}+W^{r}_{k,l}
\end{equation}
for all $\alpha\in\mathbb{N}^n$ such that
$\alpha_{r+1}-\alpha_r=k-l$, which implies
\begin{equation}
{\cal A}^{r}_{k',k'',0,l}\subset
U(P(n-1)_{\bar{0}})(f)+V^{r}_{k,l}+W^{r}_{k,l}\qquad \mbox{ for
}k''-k'+l=k.
\end{equation}
Again by induction on $t$, we obtain
\begin{equation}
{\cal A}^{r}_{k',k'',t,s}\subset
U(P(n-1)_{\bar{0}}).f+V^{r}_{k,l}+W^{r}_{k,l}\;\mbox{ for
}k''-k'+l=k,s+t=l,t\geq 0.
\end{equation}
Anyway, we get
\begin{equation}
{\cal
A}^{r}_{k,l}=U(P(n-1)_{\bar{0}})(f)+(V^{r}_{k,l}+W^{r}_{k,l}),
\end{equation}
that is $U(P(n-1)_{\bar{0}})(\bar{f})={\cal
A}^{r}_{k,l}/(V^{r}_{k,l}+W^{r}_{k,l})$.

It can be similarly proved that
$(V^{r}_{k,l}+W^{r}_{k,l})/V^{r}_{k,l}$,
$V^{r}_{k,l}/V^{r}_{k,l}\bigcap W^{r}_{k,l}$ and
$V^{r}_{k,l}\bigcap W^{r}_{k,l}$ are irreducible
$P(n-1)_{\bar{0}}$-modules.
\end{proof}
\psp

Denote
\begin{equation}
V_k^{r}=\left\{
\begin{array}{lll}
U(P(n-1))(x_{r+1}^k) & \mbox{if} & k>0,\\
U(P(n-1))(x_r^{-k}) & \mbox{if} & k\leq 0.
\end{array}
\right.
\end{equation}

\begin{thm} The module
${\cal A}_k^{r}$ has the following composition series:
\begin{equation}
{\cal A}_k^{r}\supset V_k^{r}\supset\{0\} \mbox{ if }k\neq n-r;
\end{equation}
\begin{equation}
{\cal A}_k^{r}\supset\la \theta_{r+1}\cdots\theta_n\ra\supset
V_k^{r}\supset\{0\}
\end{equation}
\end{thm}
\begin{proof}
Suppose $k\geq 0$. We prove the theorem step by step.

\textsl{1) $V^{r}_k$ is the minimal submodule of ${\cal
A}^{r}_k$.}

\psp  Let $0\neq
f(x_1,\cdots,x_n;\theta_1,\cdots,\theta_n)\in{\cal A}^{r}_k$.
Applying $E_{i,n+i}\mid_{\cal A}=x_i\partial_{\theta_i}$ for
$r<i\leq n$, and
\begin{equation}(E_{i,n+j}+E_{j,n+i})\mid_{\cal
A}=\partial_{x_i}\partial_{\theta_j}+x_j\partial_{\theta_i}\ \
\mbox{for}\; 1\leq i\leq r<j\leq n,
\end{equation}
we can get some $0\neq f_1(x_1,\cdots,x_n)\in \la f\ra$. Using
\begin{equation}
(E_{i,r}-E_{n+r,n+i})(f_1)=-x_r\partial_{x_i}(f_1)\qquad\mbox{for}\;
i\in\overline{1,r-1}\end{equation} and
\begin{equation}(E_{r+1,j}-E_{n+j,n+r+1})(f_1)=x_{r+1}\partial_{x_j}(f_1)
\qquad\mbox{for}\; j\in\overline{r+2,n},\end{equation}
 we get some $0\neq f_2(x_r,x_{r+1})\in\la f\ra$. Since
\begin{equation}
(E_{r,r+1}-E_{n+r+1,n+r})(f_2)=\partial_{x_r}\partial_{x_{r+1}}(f_2),
\end{equation}
we obtain $x_r^{\alpha_r}x_{r+1}^{\alpha_{r+1}}\in \la f\ra$ with
$\alpha_{r+1}\alpha_r=0$ and  $\alpha_{r+1}-\alpha_r=k$.\psp

\textsl{2) $V_k^{r}=\bigoplus\limits_{l=0}^{n-1}V^{r}_{k,l}$.}

\psp Assume $k\geq n-r>0$. Since
\begin{equation}
x_{r+1}^{k-l}\theta_{n-l+1}\cdots\theta_n=\frac{(k-l)!}{k!}\prod\limits_{j=n-l+1}^n(E_{n+j,r+1}-E_{n+r+1,j})(x_{r+1}^k)\in
V^{r}_k\end{equation} for $0<l<n-r$, we get $V^{r}_{k,l}\subset
V^{r}_k$ for $0<l<n-r$. Note
\begin{eqnarray}
&&(E_{n+r,r+1}-E_{n+r+1,r})(x_{r+1}^{k_n+r+1}\theta_{r+2}\cdots\theta_n)\nonumber\\
&=&(k-n+r+1)\theta_rx_{r+1}^{k-n+r}\theta_{r+2}\cdots\theta_n+x_rx_{r+1}^{k-n+r+1}\theta_{r+1}\cdots\theta_n\in
V^r_k
\end{eqnarray}
and
\begin{eqnarray}
&&(k-l+1)x_{r+1}^{k-l}\theta_{n-l}\cdots\theta_r\theta_{r+2}\cdots\theta_n\nonumber\\ & &+\sum\limits_{p=n-l}^r(-1)^{p-r}x_p\theta_{n-l}\cdots\hat{\theta_p}\cdots\theta_rx_{r+1}^{k-l+1}\theta_{r+1}\cdots\theta_n\nonumber\\
=&&\frac{1}{k-l+2}(E_{2n-l,r+1}-E_{n+r+1,n-l})\big((k-l+2)x_{r+1}^{k-l+1}\theta_{n-l+1}\cdots\theta_r\theta_{r+2}\cdots\theta_n\nonumber\\
&&+\sum\limits_{p=n-l+1}^r(-1)^{p-r}x_p\theta_{n-l+1}\cdots\hat{\theta_p}\cdots
\theta_rx_{r+1}^{k-l+2}\theta_{r+1}\cdots\theta_n\big).
\end{eqnarray}
Thus we have $V^{r}_{k,l}\subset V^{r}_k$ for $n-r\leq l\leq k+1$
such that $l\neq \frac{1}{2}(k+n-r)$. Moreover,
\begin{eqnarray}
&&\prod\limits_{j=n-l}^{n-k-2}(E_{n+r,j}-E_{n+j,r})\big(\sum\limits_{p=n-k-1}^r(-1)^{p-r}x_p\theta_{n-k-1}\cdots\hat{\theta_p}\cdots\theta_r\theta_{r+1}\cdots\theta_n\big)
\nonumber\\
&=&\sum\limits_{p=n-l}^r(-1)^{p-r}x_p\theta_{n-l}\cdots\hat{\theta_p}\cdots\theta_r\theta_{r+1}\cdots\theta_n,
\end{eqnarray}
which implies $V^{r}_{k,l}\subset V^{r}_k$ for $k+1<l<n$.

When $k+n-r$ is even, we set $l'=\frac{1}{2}(k+n-r)$. Since
$V^r_{k,l'+1}\subset V^r_k$, we have
\begin{eqnarray}& &
(k-l')x_{r+1}^{k-l'-1}
\theta_{n-l'-1}\cdots\theta_r\theta_{r+2}\cdots\theta_n \nonumber\\
&
&+\sum\limits_{p=n-l'-1}^r(-1)^{p-r}x_p\theta_{n-l'-1}\cdots\hat{\theta_p}\cdots\theta_rx_{r+1}^{k-l'}\theta_{r+1}\cdots\theta_n
\in V^r_k.
\end{eqnarray}
Note
\begin{eqnarray}
&&E_{r+1,n+r+1}\big[(k-l')x_{r+1}^{k-l'-1}
\theta_{n-l'-1}\cdots\theta_r\theta_{r+2}\cdots\theta_n\nonumber\\
&&+\sum\limits_{p=n-l'-1}^r(-1)^{p-r}x_p\theta_{n-l'-1}\cdots\hat{\theta_p}\cdots\theta_rx_{r+1}^{k-l'}\theta_{r+1}\cdots\theta_n\big]\nonumber\\
&=&(-1)^{l'-n+r}\sum\limits_{p=n-l'-1}^r(-1)^{p-r}
x_p\theta_{n-l'-1}\cdots\hat{\theta_p}\cdots\theta_rx_{r+1}^{k-l'}\theta_{r+2}\cdots\theta_n.
\end{eqnarray}
Thus $V^{r}_{k,l'}\subset V^{r}_k$ by Table 1 and (2.113).

It can be verified that
$P(n-1)\big[\sum\limits_{l=0}^{n-1}V^{r}_{k,l}\big]\subset\sum\limits_{
l=0}^{n-1}V^{r}_{k,l}$. Therefore, $V^r_k=\bigoplus\limits_{
l=0}^{n-1}V^{r}_{k,l}$.

It can be similarly proved when $k<n-r$.\psp

\textsl{3) ${\cal A}^{r}_k/V^{r}_k$ is irreducible when $k\neq
n-r$.}

\psp Again we assume $k>n-r$. The proof for $k<n-r$ is quite
similar. Let $0\neq f\in{\cal A}^{r}_{k}\setminus V^{r}_k$. We can
write
\begin{equation}
f=\sum_{k'+k''=k,t\in\ol{0,r},s\in\ol{0,n-r}}f_{k',k'',s,t},\end{equation}
where $f_{k',k'',s,t}\in{\cal A}^r_{k',k'',s,t}$ (cf.
(\ref{p7161}) and (\ref{p7162})) and only finite
$f_{k',k'',s,t}\neq0$. Since ${\cal A}^r_{k',k'',s,t}$ are all
finite dimensional $(L_1+L_2)$-modules,  $U(L_1+L_2)(f)$ is finite
dimensional. Thus there should be a weight vector $f'\in
U(L_1+L_2)(f)\subset\la f\ra$ such that $L_1^+(f')=0$ and
$L_2^+(f')=0$. Up to a scalar multiple, $f'$ should be in form of
\begin{equation}
x_r^{t-k'}\theta_{r-t+1}\cdots\theta_r
x_{r+1}^{k''-s-1}\xi_r\theta_{n-s+2}\cdots\theta_n\end{equation}
or
\begin{equation}
x_r^{t-k'}\theta_{r-t+1}\cdots\theta_rx_{r+1}^{k''-s}\theta_{n-s+1}
\cdots\theta_n
\end{equation}
for some $t\in\ol{0,r}$, $s\in\ol{0,n-r-1}$ and
$k',k''\in\mathbb{N}$ such that $k'+k''=k$. Since
\begin{eqnarray}&
&E_{r+1,n+r+1}(x_r^{t-k'}\theta_{r-t+1}\cdots\theta_r
x_{r+1}^{k''-s-1}\xi_r\theta_{n-s+2}\cdots\theta_n)\nonumber\\
&=&(-1)^tx_r^{t-k'}\theta_{r-t+1}\cdots\theta_rx_{r+1}^{k''-s}\theta_{n-s+1}
\cdots\theta_n,
\end{eqnarray}
we can assume
$f'=x_r^{t-k'}\theta_{r-t+1}\cdots\theta_rx_{r+1}^{k''-s}\theta_{n-s+1}
\cdots\theta_n$. Let $l_0=s+t$. We divide our argument into three
subcases.

(a) If $n-r\leq l_0\leq k$,
 then we have
\begin{eqnarray}
g&=&x_{r+1}^{k-l_0}\theta_{n-l_0+1}\cdots\theta_n\nonumber\\
&\equiv&\frac{(k-l_0)!}{(t-k')!(k''-s)!}(E_{r,r+1}-E_{n+r+1,n+r})^{t-k'}\prod_{i=1}^{n-s-r}(-1)^t(E_{r-t+i,r+i}\nonumber\\
&&-E_{n+r+i,n+r-t+i})(f')\;(\mbox{mod}\;V^r_k)
\end{eqnarray}
(cf. (\ref{p7165})). Note
\begin{equation}
\prod\limits_{i=n-l+1}^{n-l_0}(E_{n+i,r+1}-E_{n+r+1,i})(g)=\frac{(k-l_0)!}{(k-l)!}x_{r+1}^{k-l}\theta_{n-l+1}\cdots\theta_n\in
\la f\ra\end{equation} for $l_0<l\leq k$. Moreover,
\begin{equation}
\prod\limits_{i=n-l+1}^{n-k}(E_{n+r,i}-E_{n+i,r})(\theta_{n-k+1}\cdots\theta_n))=x_r^{l-k}\theta_{n-l+1}\cdots\theta_n\in
\la f\ra\end{equation} for $k<l\leq n$. Furthermore,
\begin{equation}
\prod\limits_{i=l}^{l_0}(E_{r+1,2n-i+1}+E_{n-i+1,n+r+1})(g)=x_{r+1}^{k-l}\theta_{n-l+1}\cdots\theta_n\in
\la f\ra\end{equation}
 for $n-r\leq l<l_0$. Thus we have
\begin{equation}x_{r+1}\xi_r\theta_{n-l+2}\cdots\theta_n\in \la
x_{r+1}^{k-n+r}\theta_{r+1}\cdots\theta_n\ra\qquad\mbox{for}\;l<n-r\end{equation}
by Theorem 2.3. Therefore, $W^{r}_{k,l}\subset \la f\ra$ for all
$0\leq l\leq n$.

(b) If $l_0\geq k$, then let
\begin{eqnarray}
g&=&x_r^{l_0-k}\theta_{n-l_0+1}\cdots\theta_n\nonumber\\
&=&\frac{(l_0-k)!}{(t-k')!(k''-s)!}(E_{r,r+1}-E_{n+r+1,n+r})^{k''-s}\prod_{i=1}^{n-s-r}(-1)^t(E_{r-t+i,r+i}\nonumber\\
&&-E_{n+r+i,n+r-t+i})(f').
\end{eqnarray}
We have
\begin{equation}\label{p7163}
\prod\limits_{i=l}^{l_0}(E_{r,2n-i+1}+E_{n-i+1,n+r})(g)=\frac{(l-l_0)!}{(l-k)!}x_r^{l-k}\theta_{n-l+1}\cdots\theta_n\in
\la f\ra\end{equation}
 for $k\leq l<l_0$. By (\ref{p7163}) and (a) with $g$ replaced by $\theta_{n-k+1}\cdots\theta_n$,
we get $W^{r}_{k,l}\subset \la f\ra$ for all $0\leq l\leq n$.

(c) If $l_0<n-r$, then
\begin{eqnarray}
&&x_{r+1}^{k-l_0}\theta_{r+1}\theta_{n-l_0+2}\cdots\theta_n\nonumber\\
&=&-\frac{(k-l_0)!}{(t-k'+1)!(k''-s)!}(E_{r,r+1}-E_{n+r+1,n+r})^{t-k'+1}\prod_{i=1}^{t-1}(-1)^t(E_{r-t+i,n-l_0+1+i}\nonumber\\
&&-E_{2n-l_0+1+i,n+r-t+i})(f')\nonumber\\
&\equiv&\frac{k-l_0}{k+n-l_0}x_{r+1}^{k-l_0-1}\xi_r\theta_{n-l_0+2}\cdots\theta_n\;(\mbox{mod}\;V^r_k)
\end{eqnarray}
(cf (\ref{p7171}). Thus
\begin{eqnarray}
g=x_{r+1}^{k-l_0-1}\xi_r\theta_{n-l_0+2}\cdots\theta_n\in\la f\ra.
\end{eqnarray}
Consequently, $x_{r+1}^{k-l-l}\xi_r\theta_{n-l+2}\cdots\theta_n\in
\la f\ra$ for all $l<n-r$ by Theorem 2.3.  Thus we obtain
$W^{r}_{k,l}\subset \la f\ra$ for $l\in\overline{0,n}$ by (a) with
$g$ replaced by $x_{r+1}^{k-n+r}\theta_{r+1}\cdots\theta_n$. \psp

Denote
\begin{equation}
W^r_k=\bigoplus\limits_{l=0}^nW^r_{k,l}.
\end{equation}
Note that $\la f\ra\supset W^{r}_k+V^{r}_k={\cal A}^{r}_k$ if
$k+n-r$ is odd. When $k+n-r$ is even, let $l'=\frac{1}{2}(k+n-r)$.
Observe that $n-r<l'<k$ and
$$f'=x_{r+1}^{k-l'-1}\theta_{n-l'}\cdots\theta_n\in
W^r_{k,l'+1}\subset\la f\ra.$$ Since
\begin{equation}
E_{r+1,n+r+1}(f')=(-1)^{r-n+l'+1}x_{r+1}^{k-l'}\theta_{n-l'}\cdots\theta_r\theta_{r+2}\cdots\theta_n,
\end{equation}
we obtain ${\cal A}^{r}_{k,l'}\subset \la f\ra$, which implies
${\cal A}^{r}_k\subset \la f\ra$.\psp

4) $k=n-r$.

Note
\begin{equation}
(E_{i,j}-E_{n+j,n+i})(\theta_{r+1}\cdots\theta_n)=0
\end{equation}
if $i,j\in\overline{1,r}$ or $i,j\in\overline{r+1,n}$ or
$i\in\overline{1,r}$ and $j\in\overline{r+1,n}$. When
$i\in\overline{r+1,n}$ and $j\in\overline{1,r}$, we have
\begin{eqnarray}
&&(E_{i,j}-E_{n+j,n+i})(\theta_{r+1}\cdots\theta_n)\nonumber\\
&=&-x_ix_j\theta_{r+1}\cdots\theta_n+(-1)^{i-r}\theta_j\theta_{r+1}\cdots\hat{\theta_i}\cdots\theta_n\nonumber\\
&=&(E_{r,j}-E_{n+j,n+r})(E_{i,r+1}-E_{n+r+1,n+i})(E_{r-1,r+1}-E_{n+r+1,n+r-1})\nonumber\\
&&\big((x_{r_1}\theta_r-x_r\theta_{r-1})x_{r+1}\theta_{r+2}\cdots\theta_n\big)\in
V^r_{n-r}.
\end{eqnarray}
Thus $P(n-1)_{\bar{0}}(\theta_{r+1}\cdots\theta_n)\subset
V^r_{n-r}$. Since
\begin{equation}
(E_{i,n+j}+E_{j,n+i})(\theta_{r+1}\cdots\theta_n)=0
\end{equation}
and
\begin{equation}
(E_{n+i,j}-E_{n+j,i})(\theta_{r+1}\cdots\theta_n)=(x_i\theta_j-x_j\theta_i)\theta_{r+1}\cdots\theta_n\in
V^r_{n-r,n-r+1},
\end{equation}
we obtain $P(n-1)_{\bar{1}}(\theta_{r+1}\cdots\theta_n)\subset
V^{r}_{n-r}$, that is,
$\la\theta_{r+1}\cdots\theta_n)\ra/V^{r}_{n-r}=\mathbb{C}\theta_{r+1}\cdots\theta_n$.
By the similar arguments as in 3), we get that ${\cal
A}^{r}_{n-r}/\la\theta_{r+1}\cdots\theta_n\ra$ is irreducible.
\end{proof}

\section{Proof of Theorem 2}
\setcounter{equation}{0} In this section, we investigate the
polynomial representation of the Lie superalgebra
$\tilde{Q}(n-1)\;(n\geq3)$.

Recall
\begin{eqnarray}
\tilde{Q}(n-1)_{\bar{0}}=\sum\limits_{1\leq i,j\leq
n}\mathbb{C}(E_{i,j}+E_{n+i,n+j}),\;\tilde{Q}(n-1)_{\bar{0}}^+=\sum\limits_{1\leq
i<j\leq n}\mathbb{C}(E_{i,j}+E_{n+i,n+j}),
\end{eqnarray}
\begin{eqnarray}
\tilde{Q}(n-1)_{\bar{1}}&=&\sum\limits_{i=1}^{n-1}\mathbb{C}(E_{i,n+i}+E_{n+i,i}-E_{i+1,n+i+1}-E_{n+i+1,i+1})\nonumber\\
&&+\sum\limits_{ i,j=1,i\neq j}^n\mathbb{C}(E_{i,n+j}+E_{n+i,j})
\end{eqnarray}
and
\begin{equation}{\cal
A}_{k,t}^r=\mbox{Span}\{x^\alpha\theta_{i_1}\cdots\theta_{i_t}\mid
i_1,\cdots,i_t\in\overline{1,n};\alpha\in\mathbb{N}^n,\;
-\sum\limits_{i=1}^r\alpha_i+\sum\limits_{j=r+1}^n\alpha_j=k-t\}.
\end{equation}
Note
\begin{equation}
H=\sum\limits_{i=1}^n\mathbb{C}E_{i,i}
\end{equation}
forms a Cartan subalgebra of $\tilde{Q}(n-1)$. We study the
representation case by case. \psp

 \textit{{Case 1}, $r=0$.}
\psp

 Set
\begin{equation}
v_t=\sum\limits_{i=1}^t(-1)^{i-1}x_i\theta_1\cdots\hat{\theta_i}\cdots\theta_t.
\end{equation}
\begin{lem}The subspace \begin{equation}{\cal A}^r_{k;t}=U(Q(n-1)_{\bar0})
(x_1^{k-t}\theta_1\cdots\theta_t)\oplus
U(Q(n-1)_{\bar0})(x_1^{k-t-1}v_{t+1})\end{equation} as
$Q(n-1)_{\bar0}$-submodules.\hfill$\square$
\end{lem}

\begin{thm} The subspace
${\cal A}^{0}_k$ has only two nonzero proper submodules $\la
x_1^k+\sqrt{k}x_1^{k-1}\theta_1\ra$ and $\la
x_1^k-\sqrt{k}x_1^{k-1}\theta_1\ra$. Moreover, $\la x_1^k\pm
\sqrt{k}x_1^{k-1}\theta_1\ra$ have a basis
\begin{eqnarray}& &
\{(k-t)x^{\alpha}\theta_{i_1}\cdots\theta_{i_t}+\sum\limits_{p=1}^t\sum\limits_{s=1}^n(-1)^p\alpha_s
x^{\alpha+\epsilon_{i_p}-\epsilon_s}\theta_s\theta_{i_1}\cdots\hat{\theta_{i_p}}\cdots\theta_{i_t}\nonumber\\
& &\pm
\sqrt{k}\sum\limits_{s=1}^n\alpha_sx^{\alpha-\epsilon_s}\theta_{i_1}\cdots\theta_{i_t}\theta_s
\mid \alpha\in\mathbb{N}^n,\;|\alpha|=k-t;\;0\leq t<\mbox{\it
min}\{k,n\};\nonumber\\ & &\;i_1,\cdots,i_t\in\overline{1,n},
\;\alpha_j>0\; \mbox{\it for  some} \;j>i_1,\cdots,i_t\}\label{q2}
\end{eqnarray}
\end{thm}
\begin{proof}
\textsl{(1) For any $0\neq f\in {\cal A}^{0}_k$, we claim that
\begin{equation}
x_1^k+\sqrt{k}x_1^{k-1}\theta_1\in\la f\ra\ \ \mbox{or}\ \
x_1^k-\sqrt{k}x_1^{k-1}\theta_1\in\la f\ra.\label{q1}
\end{equation}}

\psp Write $f=\sum_{t=0}^k (f_t+g_t)$ with
\begin{equation}f_t\in
U(Q(n-1)_{\bar{0}})(x_1^{k-t}\theta_1\cdots\theta_t),\;\;g_t\in
U(Q(n-1)_{\bar{0}})(x_1^{k-t}v_t).\end{equation} Applying
$\tilde{Q}(n-1)_{\bar{0}}^+$ to $f$, we can assume
\begin{equation}f=\sum_{t=0}^k (a_t
x_1^{k-t}\theta_1\cdots\theta_t+b_tx_1^{k-t}v_t),\qquad
a_t,b_t\in\mathbb{C}.\end{equation} Since $a_t
x_1^{k-t}\theta_1\cdots\theta_t+b_tx_1^{k-t}v_t$ and $a_{t'}
x_1^{k-t'}\theta_1\cdots\theta_{t'}+b_{t'}x_1^{k-t'}v_{t'}$ are in
different weight subspaces if $t\neq t'$, we have $a_t
x_1^{k-t}\theta_1\cdots\theta_t+b_tx_1^{k-t}v_t\in \la f\ra$. Denote
\begin{equation}t_0=\mbox{min}\{t\in\ol{0,k}\mid (a_t,b_t)\neq  (0,0)\}.\end{equation}
Observe
\begin{eqnarray}
&&(E_{1,n+2}+E_{n+1,2})\cdots(E_{1,n+t_0}+E_{n+1,t_0})(a_{t_0}x_1^{k-t_0}\theta_1\cdots\theta_{t_0}+b_{t_0}x_1^{k-t_0}v_{t_0})\nonumber\\
=&&(-1)^{\frac{t_0(t_0-1)}{2}}a_{t_0}x_1^{k-1}\theta_1+(-1)^{\frac{(t_0-1)(t_0-2)}{2}}b_{t_0}x_1^k\in\la
f\ra.
\end{eqnarray}
If $a_{t_0}\neq 0$ and $b_{t_0}=0$, we have
$x_1^{k-1}\theta_1\in\la f\ra$ and
\begin{equation}
x_1^k=(E_{1,n+1}+E_{n+1,1}-E_{2,n+2}-E_{n+2,2})(x_1^{k-1}\theta_1)\in\la
f\ra.
\end{equation}
When $a_{t_0}=0$ and $b_{t_0}\neq0$, we get $x_1^k\in\la f\ra$ and
\begin{equation}
x_1^k\theta_1=\frac{1}{k}(E_{1,n+1}+E_{n+1,1}-E_{2,n+2}-E_{n+2,2})(x_1^{k})\in\la
f\ra.
\end{equation}
Thus, under the above assumptions, (3.8) holds. In the case
$a_{t_0}\neq0$ and $b_{t_0}\neq0$, we have
$a'x_1^{k_1}\theta_1+b'x_1^k\in\la f\ra$ with
$a'=(-1)^{\frac{t_0(t_0-1)}{2}}a_{t_0}\neq0$ and
$b'=(-1)^{\frac{(t_0-1)(t_0-2)}{2}}b_{t_0}\neq0$. Since
\begin{equation}
(E_{1,n+1}+E_{n+1,1}-E_{2,n+2}-E_{n+2,2})(a'x_1^{k-1}\theta_1+b'x_1^k)=kb'x_1^{k-1}\theta_1+a'x_1^k\in\la
f\ra,
\end{equation}
we obtain $x_1^k,x_1^{k-1}\theta_1\in\la f\ra$ if
$\frac{a'}{kb'}\neq\frac{b'}{a'}$, which implies (3.8). When
$\frac{a'}{kb'}=\frac{b'}{a'}$, we have $a'/b'=\pm\sqrt{k}$, and
so (\ref{q1}) holds.\psp

\textsl{(2) The set (\ref{q2}) is a subset of $\la
x_1^k\pm\sqrt{k}x_k^{k-1}\theta_1\ra$.}

\psp Denote
\begin{equation}
\begin{array}{rl}
h(\alpha;i_1,\cdots,i_t)=&(k-t)x^{\alpha}\theta_{i_1}\cdots\theta_{i_t}+\sum\limits_{p=1}^t\sum\limits_{s=1}^n(-1)^p\alpha_s
x^{\alpha+\epsilon_{i_p}-\epsilon_s}
\theta_s\theta_{i_1}\cdots\hat{\theta_{i_l}}\cdots\theta_{i_t}\\
&+\sqrt{k}\sum\limits_{s=1}^n\alpha_sx^{\alpha-\epsilon_{i_1}-\cdots-\epsilon_{i_t}-\epsilon_s}\theta_{i_1}\cdots\theta_{i_t}\theta_s,
\end{array}
\end{equation}
where $t\in\overline{0,n-1}$. We write $h(\alpha)$ instead of
$h(\alpha;i_1\cdots,i_t)$ when $t=0$.

Since
\begin{eqnarray}
h(\alpha)=kx^\alpha+\sqrt{k}\sum\limits_{i=1}^n\alpha_ix^{\alpha-\epsilon_i}\theta_i
=\frac{\alpha_1!}{(k-1)!}\prod\limits_{i=2}^n(E_{i,1}+E_{n+i,n+1})^{\alpha_i}(x_1^k+\sqrt{k}x_1^{k-1}\theta_1)
\end{eqnarray}
we have
\begin{equation}h(\alpha)\in\la
x_1^k+\sqrt{k}x_1^{k-1}\theta_1\ra\qquad\mbox{for}\;\alpha\in\mathbb{N}^n\;\mbox{such
that}\; |\alpha|=k.\end{equation} Now we assume
\begin{equation}h(\alpha;i_1,\cdots,i_r)\in\la x_1^k+\sqrt{k}x_1^{k-1}\theta_1\ra\end{equation}
for $r<t$, $\alpha\in\mathbb{N}^n$ with $|\alpha|=k-r$, and
$i_1,\cdots,i_r\in\overline{1,n}$. Since $t<n$, we can take
$\ol{1,n}\ni j\neq i_1,\cdots,i_t$.

\begin{eqnarray}& &
h(\alpha;i_1,\cdots,i_t)\nonumber\\&=&\frac{1}{\alpha_j+1}(E_{i_1,n+j}+E_{n+i_1,j})[
h(\alpha-\epsilon_{i_1}+\epsilon_j;i_2,\cdots,i_t)]\nonumber\\
&&+\frac{1}{k}\sum\limits_{l=2}^t(-1)^lh(\alpha;i_2,\cdots,\hat{i_l},\cdots,i_t)
+\frac{(-1)^t}{\sqrt{k}}h(\alpha;i_2,\cdots,i_t).
\end{eqnarray}
So $h(\alpha;i_1,\cdots,i_t)\in\la
x_1^k+\sqrt{k}x_1^{k-1}\theta_1\ra$. By induction on $t$, the
conclusion holds.\psp

\textsl{(3) The set (\ref{q2}) forms a basis $\la x_1^k\pm
\sqrt{k}x_1^{k-1}\theta_1\ra$.}

\psp Denote by $V$ the subspace spanned by (\ref{q2}). Since
\begin{equation}
h(\alpha;i_1,\cdots,i_t)=\frac{(-1)^t}{\alpha_{i_t}+1}\sum\limits_{s\neq
i_1,\cdots,i_t}\alpha_sh(\alpha+\epsilon_{i_t}-\epsilon_s;s,i_1,\cdots,i_{t-1})
\end{equation}
if $i_t=\mbox{max}\{j\in\ol{1,n}\mid \alpha_j>0\}$, we obtain
\begin{equation}
V=\mbox{Span}\{h(\alpha;i_1,\cdots,i_t)\mid 0\leq
t<k,n;i_1,\cdots,i_t\in\overline{1,n};\alpha\in\mathbb{N}^n,\;|\alpha|=k-t\}
\end{equation}
For $j\notin\{i_1,\cdots,i_t\}$, we have
\begin{equation}
(E_{i,j}+E_{n+i,n+j})[h(\alpha;i_1,\cdots,i_t)]=\alpha_jh(\alpha+\epsilon_i-\epsilon_j;i_1,\cdots,i_t).
\end{equation}
For $i\notin\{i_1,\cdots,i_t\}$ and $j\in\{i_1,\cdots,i_t\}$, we can
assume $j=i_1$ and get
\begin{eqnarray}& &(E_{i,j}+E_{n+i,n+j})[h(\alpha;i_1,\cdots,i_t)]\nonumber\\&=&h(\alpha+\epsilon_i-\epsilon_{i_1};i,i_2,\cdots,i_t)
+\alpha_{i_1}h(\alpha+\epsilon_i-\epsilon_{i_1};i_1,\cdots,i_t).\end{eqnarray}
When $i,j\in\{i_1,\cdots,i_t\}$, we may assume $j=i_1,i=i_2$ and
have
\begin{equation}
(E_{i,j}+E_{n+i,n+j})[h(\alpha;i_1,\cdots,i_t)]
=\alpha_{i_1}h(\alpha-\epsilon_{i_1}+\epsilon_{i_2};i_1,\cdots,i_t).\end{equation}
Therefore, $(E_{i,j}+E_{n+i,n+j})(V)\subset V$ for any
$i,j\in\ol{1,n}$. Observe that
\begin{eqnarray}& &
(E_{i,n+j}+E_{n+i,j})[h(\alpha;i_1,\cdots,i_t)]\nonumber\\&=&\alpha_j\big(h(\alpha-\epsilon_j+\epsilon_i;i,i_1,\cdots,i_t)\nonumber\\
&&+\frac{1}{k}\sum\limits_{l=1}^t(-1)^lh(\alpha-\epsilon_j+\epsilon_i;i_1,\cdots,\hat{i_l},\cdots,i_t)\nonumber\\
&&+\frac{(-1)^t}{\sqrt{k}}h(\alpha-\epsilon_j+\epsilon_i;i_1,\cdots,i_t)\big)
\end{eqnarray}
if $j\notin\{i_1,\cdots,i_t\}$. When $j\in\{i_1,\cdots,i_t\}$, we
may assume $j=i_1$ and have
\begin{eqnarray}&
&(E_{i,n+i_1}+E_{n+i,i_1})[h(\alpha;i_1,\cdots,i_t)]\nonumber
\\ &=&-(E_{i,n+i_1}+E_{n+i,i_1})\big[\sum\limits_{s\neq
i_1,\cdots,i_t}h(\alpha;s,i_2,\cdots,i_t)\big].\end{eqnarray} Thus
the subspace $V$ is a submodule of $\la
x_1^k+\sqrt{k}x_1^{k-1}\theta_1\ra$, which implies $V=\la
x_1^k+\sqrt{k}x_1^{k-1}\theta_1\ra$.

Using Lemma 2.1, we get the linear independence of (\ref{q2}).
\end{proof}\psp

\psp \textit{Case 2. $1\leq r<n$.}\psp

In this case, we have:\psp

\begin{thm}
The submodule ${\cal A}_k^{r}$ is irreducible if $0<r<n$.
\end{thm}
\begin{proof}
\textsl{(1) First we claim that
\begin{equation}
\la x_r^{r+i_0}\theta_1\cdots\theta_rx_{r+1}^{k+i_0}\ra={\cal
A}_k^{r}\qquad\mbox{ for all }i_0\geq -k,-r.
\end{equation}}

\psp Since
\begin{equation}
(E_{r,r+1}+E_{n+r,n+r+1})|_{{\cal
A}^r}=\partial_{x_r}\partial_{x_{r+1}},
\end{equation}
we get
\begin{equation}
x_r^{r+i}\theta_1\cdots\theta_rx_{r+1}^{k+i}\in\la
x_r^{r+i_0}\theta_1\cdots\theta_rx_{r+1}^{k+i_0}\ra\qquad\mbox{for
}i\leq i_0.
\end{equation}
By Theorem 3.2, we derive
\begin{equation}
x^\alpha\theta_{i_1}\cdots\theta_{i_t}\theta_{j_1}\cdots\theta_{j_s}\in\la
x_r^{r+i}\theta_1\cdots\theta_rx_{r+1}^{k+i}\ra
\end{equation}
for $0\leq t\leq r$, $0\leq s\leq n-r$, and $\alpha\in\mathbb{N}^n$
such that  $\sum\limits_{p=r+1}^n\alpha_p+s=k+i$ and
$t-\sum\limits_{q=1}^r\alpha_q=-i$. Hence
\begin{equation}
\begin{array}{l}
x^\alpha\theta_{i_1}\cdots\theta_{i_t}\theta_{j_1}\cdots\theta_{j_s}\in\la
x_r^{r+i_0}\theta_1\cdots\theta_rx_{r+1}^{k+i_0}\ra.
\end{array}
\end{equation}
Note
\begin{eqnarray}
&&(E_{r+1,r}+E_{n+r+1,n+r})(x_r^{r+j}\theta_1\cdots\theta_rx_{r+1}^{k+j})\nonumber\\
=&&-x_r^{r+j+1}\theta_1\cdots\theta_rx_{r+1}^{k+j+1}+x_r^{r+j}\theta_1\cdots\theta_{r-1}x_{r+1}^{k+j}\theta_{r+1}.
\end{eqnarray}
So we get
\begin{equation}
x_r^{r+j}\theta_1\cdots\theta_rx_{r+1}^{k+j}\in\la
x_r^{r+i_0}\theta_1\cdots\theta_rx_{r+1}^{k+i_0}\ra\qquad\mbox{for
}j>i_0,
\end{equation}
which implies
\begin{equation}
\begin{array}{l}
x^\alpha\theta_{i_1}\cdots\theta_{i_t}\theta_{j_1}\cdots\theta_{j_s}\in\la
x_r^{r+i_0}\theta_1\cdots\theta_rx_{r+1}^{k+i_0}\ra
\end{array}
\end{equation}
for all $s,t\in\ol{0,r}$ and $\alpha\in\mathbb{N}^n$ such that
$\sum\limits_{p=r+1}^n\alpha_p-\sum\limits_{q=1}^r\alpha_q=k-s-t$,
that is,  $\la
x_r^{r+i_0}\theta_1\cdots\theta_rx_{r+1}^{k+i_0}\ra={\cal
A}_k^{r}$.\psp

\textsl{(2) Next for any $0\neq f\in{\cal A}_k^{r}$, we have
\begin{equation}
x_r^{r+i}\theta_1\cdots\theta_rx_{r+1}^{k+i}\in\la f\ra \qquad\mbox{
for some }i\geq -k,-r.
\end{equation}}

In fact, we can rewrite $f=\sum\limits_i g_ih_i$, with
\begin{eqnarray}&
g_i\in \mbox{Span}\!\!\!&\{x^\beta\theta_{j_1}\cdots\theta_{j_s}\mid
0\leq s\leq n-r;\beta\in\mathbb{N}^n,\;|\beta|=k+i-s,\nonumber\\
& &\beta_1=\cdots=\beta_r=0;j_1,\cdots,j_s\in\overline{r+1,n}\},
\end{eqnarray}
\begin{eqnarray}&
h_i\in
\mbox{Span}\!\!\!&\{x^\alpha\theta_{i_1}\cdots\theta_{i_t}\mid 0\leq
t\leq r;\alpha\in\mathbb{N}^n,\;|\alpha|=i+t,\nonumber\\
& &\alpha_{r+1}=\cdots=\alpha_n=0;i_1\cdots,i_t\in\overline{1,r}\}.
\end{eqnarray}
Applying $E_{p,q}+E_{n+p,n+q}$ to $f$, we may assume
\begin{equation}
(E_{p,q}+E_{n+p,n+q})(g_i)=0\qquad\mbox{for }r<p<q\leq
n,\end{equation}
\begin{equation}(E_{p,q}+E_{n+p,n+q})(h_i)=0\qquad\mbox{for }0<p<q\leq r.
\end{equation}
Therefore, there should be some
\begin{eqnarray}
f_1&=&a(x_r^{t+i}\theta_1\cdots\theta_t)(x_{r+1}^{k+i-s}\sum\limits_{p=1}^s(-1)^px_{r+p}\theta_{r+1}\cdots
\hat{\theta}_{r+p}\cdots\theta_{r+s})\nonumber\\
&&+b(x_r^{t+i}\theta_1\cdots\theta_t\xi'_r)(x_{r+1}^{k+i-s}\sum\limits_{p=1}^s(-1)^px_{r+p}\theta_{r+1}
\cdots\hat{\theta}_{r+p}\cdots\theta_{r+s})\nonumber\\
&&+c(x_r^{t+i}\theta_1\cdots\theta_t)(x_{r+1}^{k+i-s}\theta_{r+1}\cdots\theta_{r+s})\nonumber\\
&&+d(x_r^{t+i}\theta_1\cdots\theta_t\xi'_r)(x_{r+1}^{k+i-s}\theta_{r+1}\cdots\theta_{r+s})\in\la
f\ra
\end{eqnarray}
for some $i,t,s\in\ol{1,n}$ and $a,b,c,d\in\mathbb{C}$, where
$\xi'_r=\sum\limits_{j=1}^rx_j\theta_j$. Note that if $t<r-1$, then
\begin{equation}
0\neq\prod\limits_{p=t+1}^{r-1}(-1)^p(E_{n+p,r}+E_{p,n+r})(f_1)\in\la
f\ra.\end{equation} So we can assume $t=r-1$. If $s>1$, we have
\begin{equation}
0\neq\prod\limits_{p=2}^s(E_{r+1,n+r+p}+E_{n+r+1,r+p})(f_1)\in\la
f\ra.\end{equation} Thus we can assume $s=1$. Under the assumptions,
\begin{eqnarray}
f_1&=&ax_r^{r+i-1}\theta_1\cdots\theta_{r-1}x_{r+1}^{k+i}+bx_r^{r+i}\theta_1\cdots\theta_rx_{r+1}^{k+i}\nonumber\\
&&+c(x_r^{r+i-1}\theta_1\cdots\theta_{r-1})(x_{r+1}^{k+i-1}\theta_{r+1})
+d(x_r^{r+i}\theta_1\cdots\theta_r)(x_{r+1}^{k+i-1}\theta_{r+1})
\end{eqnarray}
and
\begin{equation}
(E_{p,n+p}+E_{n+p,p})(f_1)=0\qquad\mbox{ for
}p\in\overline{1,n},p\neq r,r+1.\end{equation}

Set
\begin{eqnarray}
f_2&=&(-1)^{r-1}(E_{r+1,n+r+1}+E_{n+r+1,r+1}-E_{p,n+p}-E_{n+p,p})(f_1)\nonumber\\
&=&cx_r^{r+i-1}\theta_1\cdots\theta_{r-1}x_{r+1}^{k+i}-dx_r^{r+i}\theta_1\cdots\theta_rx_{r+1}^{k+i}\nonumber\\
&&+a(k+i)(x_r^{r+i-1}\theta_1\cdots\theta_{r-1})(x_{r+1}^{k+i-1}\theta_{r+1})
\nonumber\\ &
&-b(k+i)(x_r^{r+i}\theta_1\cdots\theta_r)(x_{r+1}^{k+i-1}\theta_{r+1}),\end{eqnarray}
\begin{eqnarray}
f_3&=&(-1)^{r-1}(E_{r,n+r}+E_{n+r,r}-E_{p,n+p}-E_{n+p,p})(f_1)\nonumber\\
&=&b(r+i)x_r^{r+i-1}\theta_1\cdots\theta_{r-1}x_{r+1}^{k+i}-ax_r^{r+i}\theta_1\cdots\theta_rx_{r+1}^{k+i}\nonumber\\
&&+d(r+i)(x_r^{r+i-1}\theta_1\cdots\theta_{r-1})(x_{r+1}^{k+i-1}\theta_{r+1})
\nonumber\\ &
&-c(x_r^{r+i}\theta_1\cdots\theta_r)(x_{r+1}^{k+i-1}\theta_{r+1}),\end{eqnarray}
\begin{eqnarray}
f_4&=&(-1)^{r-1}(E_{r+1,n+r+1}+E_{n+r+1,r+1}-E_{p,n+p}-E_{n+p,p})(f_3)\nonumber\\
&=&d(r+i)x_r^{r+i-1}\theta_1\cdots\theta_{r-1}x_{r+1}^{k+i}+cx_r^{r+i}\theta_1\cdots\theta_rx_{r+1}^{k+i}\nonumber\\
&&+b(r+i)(k+i)(x_r^{r+i-1}\theta_1\cdots\theta_{r-1})(x_{r+1}^{k+i-1}\theta_{r+1})
\nonumber\\&&+a(k+i)(x_r^{r+i}\theta_1\cdots\theta_r)(x_{r+1}^{k+i-1}\theta_{r+1}).
\end{eqnarray}

We can assume
\begin{equation}
\big(c^2-a^2(k+i)\big)+(r+i)\big(d^2-b^2(k+i)\big)\neq0\label{q4}
\end{equation}
because we can take
\begin{eqnarray}
&&-(E_{r+1,r}+E_{n+r+1,n+r})(f_1)\nonumber\\
&=&ax_r^{r+i}\theta_1\cdots\theta_{r-1}x_{r+1}^{k+i+1}+bx_r^{r+i+1}\theta_1\cdots\theta_rx_{r+1}^{k+i+1}\nonumber\\
&&+(c-b)(x_r^{r+i}\theta_1\cdots\theta_{r-1})(x_{r+1}^{k+i}\theta_{r+1})
+d(x_r^{r+i+1}\theta_1\cdots\theta_r)(x_{r+1}^{k+i}\theta_{r+1})
\end{eqnarray}
instead of $f_1$ if (\ref{q4}) does not hold.

By (3.41) and (3.46)-3.48), $f_1,f_2,f_3,f_4\in\la f\ra$. Hence
\begin{eqnarray}
f_5&=&cf_1-af_2+df_3-bf_4\nonumber\\
&=&\big(c^2-a^2(k+i)+d^2(r+i)-b^2(k+i)(r+i)\big)\nonumber\\
&&\times(x_r^{r+i-1}\theta_1\cdots\theta_{r-1})(x_{r+1}^{k+i-1}\theta_{r+1})\in\la
f\ra.
\end{eqnarray}
Thus
$(x_r^{r+i-1}\theta_1\cdots\theta_{r})(x_{r+1}^{k+i-1}\theta_{r+1})\in\la
f\ra$. Note
\begin{eqnarray}& &
x_r^{r+i}\theta_1\cdots\theta_{r-1}x_{r+1}^{k+i}\nonumber\\ &=&
(E_{r,n+r}+E_{n+r,r}-E_{p,n+p}-E_{n+p,p})(E_{r+1,n+r+1}+E_{n+r+1,r+1}\nonumber\\
&&-E_{p,n+p}-E_{n+p,p})(x_r^{r+i-1}\theta_1\cdots\theta_{r-1})(x_{r+1}^{k+i-1}\theta_{r+1}).
\end{eqnarray}
So $x_r^{r+i}\theta_1\cdots\theta_{r}x_{r+1}^{k+i}\in\la f\ra$. This
completes the proof of the theorem.
\end{proof}
\vspace{1cm}

\begin{center}{\Large \bf Acknowledgement}\end{center}

I would like to thank Professor Xiaoping Xu for his advice and
suggesting this research topic.

\vspace{1cm}

\bibliographystyle{amsplain}

\begin{thebibliography}{10}

\bibitem{} J. Brundan, Kazhdan-Lusztig polynomials and character
formulae for the Lie superalgebra $q(n)$, {\it Adv. Math.} {\bf
182} (2004), 28-77.


\bibitem{} L. Frappat, A. Sciarrino and P. Sorba, Dynkin-like diagrams
and representation of the strange superalgebra $P(n)$, {\it J. Math.
Phys.} {\bf 32} (1991), 3268-3277.

\bibitem{} L. Frappat, A. Sciarrino and P. Sorba, Oscillator
realization of the strange superalgebras $P(n)$ and $Q(n)$, {\it J.
Math. Phys.} {\bf 33} (1992), 3911

\bibitem{}\label{r4} L. Frappat, A. Sciarrino and P. Sorba, {\it Dictionary on Lie
Algebras and Superalgebras}, Academic Press, 2000

\bibitem{} M. Gorelik, The centre of simple P-type Lie superalgebra,
{\it J. Algebra} {\bf 246} (2001), 414-428.

\bibitem{} M. Gorelik, Shapovalov determinants of Q-type Lie
superalgebras, {\it IMRP} (2006), Art. Id. 96895.

\bibitem{} M. Gorelik and V. Serganova, On representations of the
affine superalgebra $Q(n)^{(2)}$, {\it Mosc. Math. J.} {\bf 8}
(2008), 91-109.

\bibitem{} C. Gruson, Sur la cohomologic des super alg\`{e}bres de
Lie \'{e}tranges, {\it Transform. Groups} {\bf 5} (2000), no. 1,
73-84.


\bibitem{} P. D. Javis and M. K. Murray, Casimir invariants,
characteristic identities, and tensor operators for ``strange
superalgebras, {\it J. Math. Phys.} {\bf 24} (1983), no. 7,
1705-1710.

\bibitem{}\label{r6} C. Luo, Noncanonical polynomial representations of
classical Lie algebras, arXiv: 0804.0305 [math.RT].

\bibitem{} M. L. Nazarov, Yangians of the ``strange" Lie
superalgebras, {\it Quantum groups} ({\it Leningrad}, 1990), 90-97,
{\it Lecture Notes in Math.,} 1510, {\it Springer, Berlin,} 1992.

\bibitem{} C. Martinez and E. Zelmanov,  Lie superalgebras
graded by $P(n)$ and $Q(n)$, {\it Proc. Natl. Acad. Sci. USA} {\bf
100} (2003), no. 14, 8130-8137.

\bibitem{} B. Medak, The group superalgebras in ``strange" Lie
superalgebras $P(n)$, {\it Hadronic J. Suppl.} {\bf 12} (1997), no.
2, 171-216.

\bibitem{} B. Medak, On the compatibility of $\mathbb{Z}$- and
$\mathbb{Z}_2$-gradations at ``strange" Lie superalgebras $P(n)$
pointed by the Jacobi identity, {\it Rep. Math. Phys.} {\bf 49}
(2002), no. 2-3, 305-314.


\bibitem{} D. Moon, Tensor product representations of the Lie
superalgebra $p(n)$ and their centralizers, {\it Commun. Algebra}
{\bf 31} (2003), no. 5, 2095-2140.

\bibitem{}V. G . Kac, Lie superalgebras, {\it Adv. Math.} {\bf 26} (1977), 8-96.

\bibitem{} V. G. Kac, Characters of typical representations of
classical Lie superalgebras, {\it Commun. Algebra} {\bf 5} (1977),
889-897.

\bibitem{} V. G. Kac, {\it Representations of Classical Lie
Superalgebras}, Lecture Notes in Math 676, Spring, Berlin, 1978,
597-626.

\bibitem{} T. Palev and J. Van der Jeugt, Fock representations of
the Lie superalgebra $q(n+1)$, {\it J. Phys. A: Math. Gen.} {\bf 33}
(2000), 2527-2544.

\bibitem{} I. Penkov and V. Serganova, Characters of finite-dimensional irreducible
$q(n)$-modules, {\it Lett. Math. Phys.} {\bf 40} (1997), no. 2,
147-158.

\bibitem{} V. Serganova, On representations of the Lie superalgebra
$p(n)$, {\it J. Algebra} {\bf 258} (2002), no. 2, 615-630.


\bibitem{} V. Stukopin, Yangian of the strange Lie superalgebra of
$Q_{n-1}$ type, Drinfel'd approach, {\it SIGMA Symmetry
Integrability Geom. Methods Appl.} {\bf 3} (2007), Paper 069, 12pp.



\end{thebibliography}

\end{document}